\documentclass[reqno,12pt]{amsart}
\usepackage[sort&compress,numbers]{natbib}
\usepackage{amsfonts}
\usepackage{mathdots,amssymb,graphicx,amssymb,amsmath,amsopn, mathtools,mathabx}
\usepackage[utf8]{inputenc}
\usepackage{enumitem}
\usepackage{tikz-cd}
\usepackage{scalerel}
\usepackage{pgfplots}
\usepackage{amscd}
\usepackage{latexsym}
\usepackage{xcolor}
\usepackage{amscd}
\usepackage{pb-diagram}
\usepackage{mathrsfs}
\usepackage{accents}
\usepackage[textwidth=15cm,textheight=22.5cm,headheight=0.6cm,headsep=1cm,centering]{geometry}
\usepackage{fancyhdr}
\usepackage{float}

\usepackage[bookmarks=true,colorlinks=true,linkcolor=airforceblue,citecolor=brightmaroon,plainpages=tfalse,pdftex]{hyperref}
\hypersetup{bookmarksopen=true,pdfview=fitbh}

\definecolor{brightmaroon}{rgb}{0.76, 0.13, 0.28}
\definecolor{airforceblue}{rgb}{0, 0.25, 0.77}
\definecolor{myOrange}{rgb}{1,0.5,0}
\definecolor{brightmaroon}{rgb}{0.76, 0.13, 0.28}
\definecolor{airforceblue}{rgb}{0, 0.4, 0.66}

\allowdisplaybreaks
\pgfplotsset{compat=1.18}

\theoremstyle{plain}

\newtheorem{teo}{Theorem}[section]

\newtheorem{pro}[teo]{Proposition}
\newtheorem{coro}[teo]{Corollary}
\newtheorem{defi}[teo]{Definition}
\newtheorem{remark}[teo]{Remark}

\numberwithin{equation}{section}
\newcommand{\Su}{\mathbf{S}}
\newcommand{\un}{\mathbf{u}}
\newcommand{\hn}{\mathbf{h}}

\newcommand{\N}{\mathbb{N}}

\numberwithin{equation}{section}
\newcommand{\prodint}[1]{\left\langle{#1}\right\rangle}

\usepackage{color}

\newcommand*\elln{\ensuremath{\boldsymbol\ell}}
\newcommand*\sn{\ensuremath{\boldsymbol s}}

\begin{document}
	
\title[Symmetrization process and truncated orthogonal polynomials]{Symmetrization process and truncated orthogonal polynomials}
	
\author[ D. Dominici, J. C García-Ardila,  Marcellán]{Diego Dominici, Juan C. García-Ardila,  Francisco Marcellán}

\address[D. Dominici]{Department of Mathematics, State University of New York at New Paltz, 1 Hawk Dr., New Paltz, NY 12561-2443, USA. 
	\\Research Institute for Symbolic Computation, Johannes Kepler University Linz, Altenberger Straße 69, 4040 Linz, Austria.}
\email{ddominic@risc.jku.at}
	
 \address[J. C. Garc\'ia-Ardila]{Departamento de Matem\'atica Aplicada a la Ingenier\'ia Industrial \\Universidad Polit\'ecnica de Madrid\\ Calle Jos\'e Gutierrez Abascal 2, 28006 Madrid, Spain.}
	\email{juancarlos.garciaa@upm.es}

	\address[F. Marcellán]{Departamento de Matemáticas, Universidad Carlos III de Madrid, Leganés, Spain} \email{pacomarc@ing.uc3m.es}
	
\thanks{ }
	
\date{\today}
\begin{abstract}
We define the family of truncated Laguerre polynomials $P_n(x;z)$, orthogonal with respect to the linear functional $\elln$ defined by
$$\prodint{\elln,p}=\int_{0}^zp(x)x^\alpha e^{-x}dx,\qquad\alpha>-1.$$

The connection between $P_n(x;z)$ and the polynomials $S_n(x;z)$ (obtained through the symmetrization process) constitutes a key element in our analysis. As a consequence, several properties of the polynomials $P_n(x;z)$ and $S_n(x;z)$ are studied taking into account the relation between the parameters of the  three-term recurrence relations that they satisfy. Asymptotic expansions of these coefficients are given.  Discrete Painlev\'e and Painlev\'e equations associated with such coefficients appear in a natural way. An electrostatic interpretation of the zeros of such polynomials as well as the dynamics of the zeros in terms of the parameter $z$ are given.
\end{abstract}
	
\maketitle
\tableofcontents
\section{Introduction}

In \cite{Chen05},  the study of the behavior of the recurrence coefficients for orthogonal polynomials with respect to singular weights was considered. There are some motivations to deal with this problem. On the one hand, the $x$ -ray problem of the condensed matter theory. On the other hand, related problems in random matrix theory which involve the asymptotics of the Fredholm determinants of finite convolution operators with discontinuous symbols. The aim is the computation for large $n$  of the determinant of the Hankel matrix $H_{n}=[\un_{j+k}]_{j, k=0}^{n-1}$ with moments $\un_{m}=\int_{a}^{b} t^{m} w(t) dt,$ $m=0,1,\ldots.$

In such a contribution, the authors considered the deformed weight $w(x)$ of a weight function $w_{0}(x)$ supported on an interval $[a,b]$ with one jump  at the point $z\in [a, b]$ in such a way that
$w(x) = w_{0}(x) (1- \frac{\beta}{2} + \beta H (x- z)),$ where $H(y)$ is the Heaviside step function and $\beta$ represents the height of the jump. For the Hermite weight $w_{0}(x)= \exp (-x^2)$ and $z=0,$ the asymptotic behavior of the coefficients of the three term recurrence relation that the corresponding sequence of orthogonal polynomials satisfy as well as some computational tests are presented. The authors stated that \emph{unfortunately, the formalism developed in this paper does not seem to shed any light on its delicate dependence on $\beta$. On the other hand,  the scheme breaks down for $z\neq 0$ and  the numerics does not suggest an ansatz for the asymptotes.}

When several  jumps are considered, a similar approach can be done. The Hermite case has been extensively studied in the literature by many authors, (see \cite{MR2265830,Lyu1,Wu}, in the framework of Gaussian unitary ensembles. On the other hand, in the particular case of  two symmetric jumps, you get the so called Rys polynomials which play an important role in the evaluation of the electron repulsion integrals in quantum chemistry computer codes (see \cite{Gradimir18,Shizgal}, among others). A connection with integrable systems in terms of a Toda perturbation of the Lebesque measure supported on $[-1,1]$ has been explored by many authors (see \cite{DM22} and the updated list of references therein).\\

In \cite{Lyu2},  the Hankel determinant generated by the Laguerre weight with jump discontinuities at $m$  points of the positive real semiaxis $z_{k},\, k = 1, \ldots , m$ wais studied. By employing the ladder operator approach to establish Riccati equations, the authors deduced that  $\sigma_{n} (z_{1}, \cdots, , z_{m})$, the logarithmic derivative of the determinant of the Hankel determinant $H_{n}$ satisfies a generalization of the $\sigma$-form of Painlev\'e V equation.  By considering the Riemann-Hilbert problem for the associated orthogonal polynomials and via the corresponding  Lax pair, $\sigma_{n}$ is given  in terms of solutions of a coupled Painlev\'e V system. The connection with the Laguerre unitary examples has been stated in  \cite{MR2525311,Basor,Lyu2}.\\

The aim of this contribution is to analyze the properties of polynomials orthogonal with respect to the truncated distribution defined  by  a jump at a positive real number $z$ and $\beta= -2.$ In other words, we deal with a linear functional $\elln$ defined in the linear space of polynomials with real coefficients such that $\prodint{\elln,p}= \int_{0}^{z} p(x) x^{\alpha} \exp (-x) dx.$  We use the fact that such a linear functional is semiclassical and, as a consequence, the powerful  tools of such a theory of linear functionals is used to get a natural approach to algebraic and analytic properties of the coefficients of the three term recurrence relation the corresponding sequences of orthogonal polynomials satisfy. Semiclassical linear functionals were introduced in \cite{Ma85,Ma87} and a comprehensive survey is \cite{Ma91}. The concept of class of a semiclassical linear functional allows to establish a hierarchy such that in the bottom the classical ones (Hermite, Laguerre, Jacobi, Bessel) appear. For the description of semiclassical linear functionals of \textit{class one}, see \cite{Belmehdi}.

For semiclassical linear functionals a system of nonlinear difference equations (Laguerre-Freud equations) that the coefficients of the three-term recurrence satisfy is deduced from the Pearson equation associated with the linear functional (see \cite{Said}). In a natural way, discrete Painlev\'e equations appear and its study has attracted the interest of many researchers, (see \cite{Magnus2,Walter}, among others). Painlev\'e differential equations for such coefficients for particular cases of semiclassical linear functionals are analyzed in \cite{Magnus1,Walter}.\\

The structure of the manuscript is as follows. In Section~\ref{section2} we present a basic background on linear functionals, orthogonal polynomials, the symmetrization process for a linear functional  and  semiclassical linear functionals.  In Section~\ref{section3}, the truncated  Laguerre linear functional is introduced. The symmetrization process yields a truncated generalized Hermite linear functional. Both of them are semiclassical and the corresponding classes are analyzed therein. In Section~\ref{section4}, the Laguerre-Freud equations, i. e., nonlinear equations satisfied by the parameters of the three term recurrence relations that the corresponding orthogonal polynomials satisfy are deduced in both cases. In Section~\ref{section5}, the ladder operators (lowering and raising operators) for the truncated Laguerre  and truncated generalized Hermite linear functionals are studied. As a consequence, the second order linear differential equation that the corresponding sequences of orthogonal polynomials satisfy are given. In Section~\ref{section6}, power series expansions of the coefficients of the three-term recurrence relations as functions of the parameter $z$ are analyzed and, in particular, an iterative method in order to find  the Taylor expansion of the coefficient $\sigma_{n}(z)$ of the monomial $ x^{n-1}$ of the $n$-th orthogonal polynomial $P_{n}(x; z)$ is given. As a consequence, in Section~\ref{section7}, the asymptotic analysis of such coefficients is deduced. In Section~\ref{section8}, discrete Painlev\'e  and Painlev\'e equations for the  parameters of the three term recurrence relation of the orthogonal polynomials  associated with the truncated Laguerre linear functional are given. In Section~\ref{section9} we focus our attention on the electrostatic interpretation of the zeros of  the above orthogonal polynomials. Their dynamical behavior is studied in terms of the parameter $z.$ As a direct consequence of the symmetrization process, the above problems are analyzed in the case of the truncated  generalized Hermite linear functional. Concluding remarks as well as some open problems are presented in Section~\ref{section10}.
\section{Basic background}\label{section2}
Let $\mathbb{N}=\{1,2,\ldots\},$ $\mathbb{N}_0=\mathbb{N}\cup\{0\}$,     $\mathbb{R}$ and $\mathbb{C}$ be the usual sets of real and complex numbers, respectively.  If $\mathbb{K}$ is a field (we mostly think of $\mathbb{K}$ as $\mathbb{R}$
	or $\mathbb{C}$), we denote the ring of  \emph{formal power series}  in the variable $z$ by
	\begin{equation*}
\mathbb{F}=\left\{
{\displaystyle\sum\limits_{n=0}^{\infty}}
c_{n}z^{n}:\quad c_{n}\in\mathbb{K}\right\}.
\end{equation*}
Let $\un$ be a complex-valued linear functional   defined on the linear space of polynomials with complex coefficients $\mathbb{P}$, that is, $$\un:\mathbb{P}\to \mathbb{C}, \quad p(x)\to \prodint{\un,p(x)}.$$  We denote the $n$th moment of $\un$ by $\un_n:=\prodint{\un,x^n}$, $n\geq0$.
\begin{defi} Let $\un$  be a linear functional. Then
\begin{enumerate}
\item For $q\in\mathbb{P}$ we define the linear functional $q(x)\un$ as 	\begin{equation*}
\prodint{q(x)\un,p(x)}=\prodint{\un,q(x)p(x)}, \quad p\in\mathbb{P},
	\end{equation*}
\item The derivative of  $\un$ is the linear functional such that
$$\prodint{D\un,p(x)}=-\prodint{\un,p^\prime(x)}.$$
\end{enumerate}
	\end{defi}
 The linear functional $\un$ is said to be quasi-definite (resp. positive-definite) if every leading principal submatrix of the Hankel matrix $H=(\un_{i+j})_{i , j=0}^{\infty}$ is nonsingular (resp. positive-definite).  In such a situation, there exists a sequence of monic polynomials $(P_n)_{n\geq0}$ such that $\deg{P_n}=n$ and $\prodint{\un,P_n(x)P_m(x)}=K_{n}\delta_{n,m},$ where $\delta_{n,m}$ is the Kronecker symbol and $K_{n}\ne 0$ (see \cite{Chi,GMM21}). The sequence $(P_n)_{n\geq0}$ is said to be the sequence of monic orthogonal polynomials (SMOP) with respect to the linear functional~$\un$. \\
	
Let $\un$ be a quasi-definite linear functional and $(P_n)_{n\geq0}$   its corresponding  SMOP. Then  there  exist two  sequences of complex numbers $( a_n)_{n\geq1}$ and $(b_n)_{n\geq0}$, with $ a_n\ne 0$, such that	
\begin{equation}\label{ttrrr}
\begin{aligned}
	x\,P_{n}(x)&=P_{n+1}(x)+b_n\,P_{n}(x)+ a_{n}\,P_{n-1}(x),\quad n\geq0,\\
P_{-1}(x)&=0, \ \ \ \ \ P_{0}(x)=1.
\end{aligned}
\end{equation}
On the contrary, by Favard's Theorem (see \cite{Chi,GMM21}) if $(P_n)_{n\geq0}$ is a sequence of monic polynomials generated by a three-term recurrence relation as in \eqref{ttrrr} with $ a_n\ne 0,$ $n\geq1$, then there exists a unique linear functional $\un$ such that $(P_n)_{n\geq0}$ is its SMOP.
\begin{defi}
Given a quasi-definite linear functional $\un$ with moments $(\un_n)_{n\geq0}$, the formal series
\begin{equation*}
\Su_{\un}(z)=:\sum_{n=0}^\infty \dfrac{\un_n}{z^{n+1}}
\end{equation*}
is said to be the \textit{Stieltjes function} associated with $\un$.
\end{defi} Due to the connection between the Stieltjes function  and the coefficients of the three term recurrence relation \eqref{ttrr} through   continued fractions, the Stieltjes function is a very useful tool to study properties of quasi-definite linear functionals \cite{MR1164865,MR2442472,MR1083352,MR0025596}.
\begin{defi}[\cite{Ma87}]
A quasi-definite functional $\mathbf{u}$ is semiclassical if there exist non-zero polynomials $\phi(x)$ and $\psi(x)$ with $\deg\phi(x)=:r\ge 0$ and $\deg\psi(x)=:t\ge 1$, such that $\mathbf{u}$ satisfies the distributional Pearson equation
\begin{equation}\label{pearson-semic}
D(\phi(x)\,\mathbf{u})+\psi(x)\,\mathbf{u}=0.\end{equation}
A sequence of orthogonal polynomials associated with $\mathbf{u}$ is called a semiclassical sequence of orthogonal polynomials.
\end{defi}
Notice that a semiclassical linear functional satisfies many Pearson equations, Indeed, if $q(x)$ is a polynomial and $\un$ satisfies \eqref{pearson-semic}, then $\un$ also satisfies the Pearson equation
$$D(\widetilde\phi(x)\,\mathbf{u})+\widetilde\psi(x)\,\mathbf{u}=0,$$
where $\widetilde\phi(x)=q(x)\phi(x)$ and $\widetilde\psi(x)=(\psi(x)q(x)-\phi(x)q^\prime(x)).$
The non-uniqueness of the Pearson equation motivates the following definition.
\begin{defi}
The class of a semiclassical functional $\mathbf{u}$ is defined as
\begin{equation*}
	\mathfrak{s}(\mathbf{u})=: \min \max\{\deg \phi(x)-2, \deg\psi(x)-1 \},
\end{equation*}
where the minimum is taken among all pairs of polynomials $\phi(x)$ and $\psi(x)$ so that \eqref{pearson-semic} holds.
\end{defi}
\begin{pro}[\cite{GMM21, Ma91}]\label{sim_cond}
Let $\mathbf{u}$ be a semi-classical linear functional and let $\phi(x)$ and $\psi(x)$ be non-zero polynomials with $\deg\phi(x)=:r$ and $\deg \psi(x)=:t$, such that \eqref{pearson-semic}  holds. Let $s =: \max(r-2,t-1)$. Then $s = \mathfrak{s}(\mathbf{u})$  if and only if
\begin{equation}\label{prodformula}	\prod_{c:\,\phi(c)=0}\left(|\psi(c)+\phi^\prime(c)|+|\langle\mathbf{u},\theta_c\psi(x)+\theta^2_c\phi(x)\rangle|\right)>0.
\end{equation}
Here, $\theta_c f(x)=\dfrac{f(x)-f(c)}{x-c}.$
\end{pro}
\begin{defi}A quasi-definite linear  functional  $\mathbf{s}$ is said to be \textit{symmetric} if
$$\prodint{\sn,x^{2k+1}}=0,\quad k\ge 0, $$
that is, all its odd moments are zero.
\end{defi} 	
The following result characterizes symmetric quasi-definite linear functionals in terms of the recurrence relation satisfied by its MOPS.
	\begin{teo}[\cite{Chi}]
		Let $(S_n(x))_{n\ge 0}$ be an MOPS associated with a quasi-definite linear functional~$\mathbf{s}$. Then the following statements are equivalent:
		\begin{enumerate}
			\item[(i)] $\mathbf{s}$ is symmetric.
			\item[(ii)] $S_{n}(x)$ has the same parity as $n$, that is, $S_{n}(x)$ is an even (resp. odd) function when $n$ is even (resp. odd).
			\item[(iii)] $(S_n(x))_{n\ge 0}$ satisfies a recurrence relation
			\begin{equation}\label{Tre}
				\begin{aligned}
					xS_n(x)&=S_{n+1}(x)+\gamma_n\,S_{n-1}(x), \ \ \ \ \ n\geq0.\\
					S_{-1}(x)&=0, \ \ \ \ \ S_{0}(x)=1.
				\end{aligned}
			\end{equation}
		\end{enumerate}
	\end{teo}
	If $(S_n(x))_{n\ge 0}$ is the MOPS associated with a linear symmetric functional, then $S_n(x)$ only contains powers of $x$ with the same parity as $n$. Therefore, there are two monic polynomial sequences $(P_n(x))_{n\ge 0}$ and $(Q_n(x))_{n\ge 0}$, such that  		
	\begin{equation}\label{pareimpar}
		S_{2n}(x)=P_n(x^2) \ \ \ \ \ \text{and}\ \ \ \ \ S_{2n+1}(x)=x\,Q_n(x^2), \quad n\geq0.
	\end{equation}
	The following result states that $(P_n(x))_{n\ge 0}$ and $(Q_n(x))_{n\ge 0}$ are also sequences of orthogonal polynomials.
	
	\begin{pro}[Chihara \cite{Chi}]
		Let $\mathbf{s}$  be a quasi-definite  symmetric linear functional. Then the sequence of monic polynomials  $(P_{n}(x))_{n\ge 0}$  is orthogonal with respect to the functional  $\un$  defined by  $\prodint{\un,x^n}=\prodint{\sn,x^{2n}}$ and the sequence of monic polynomials$(Q_{n}(x))_{n\ge 0}$ is orthogonal with respect to the functional $x\un$.
	\end{pro}
	\begin{pro}[\cite{GMM21}]
		Let $\sn$ a quasi-definite symmetric linear functional with MOPS $(S_n(x))_{n\geq 0}$ and let $(P_n(x))_{n\geq 0}$ and $(Q_n(x))_{n\geq 0}$ be the  MOPS  associated with the linear functionals  $\un$ and $x\un$, respectively, defined in \eqref{pareimpar} and satisfying the recurrence relations
		\begin{equation}\label{ttrr}
			\begin{aligned}
				x\,P_{n}(x)&=P_{n+1}(x)+b_n\,P_{n}(x)+ a_{n}\,P_{n-1}(x),\quad n \geq 0,\\
				x\,Q_{n}(x)&=Q_{n+1}(x)+d_n\,Q_{n}(x)+ c_{n}\,Q_{n-1}(x),\quad n\geq 0,
			\end{aligned}
		\end{equation}
		with $P_{-1}(x)=Q_{-1}(x)=0$, $P_{0}(x)=Q_{0}(x)=1.$ Then,
		\begin{equation*}
			\begin{aligned}
&x\,\mathbf{P}=(LU)\mathbf{P},\quad &\mathbf{P}&=(P_0,P_1,\cdots)^\top\\
&x\,\mathbf{Q}=(UL)\mathbf{Q},\quad &\mathbf{Q}&=(Q_0,Q_1,\cdots)^\top	
\end{aligned}
\end{equation*}
where    $U$  and  $L$ are upper and lower triangular matrices, respectively, given by				
		\begin{equation*}
			L=
			\begin{pmatrix}
				1       &        &  &\\
				\gamma_2&1       &  &\\
				&\gamma_4&1 &\\
				&        &\ddots&\ddots\\
			\end{pmatrix}, \ \ \ \ 	U=
			\begin{pmatrix}
				\gamma_1&1      &        & &\\
				&\gamma_3&1       & &\\
				&        &\gamma_5&1 &\\
				&        &        &\ddots&\ddots\\
			\end{pmatrix}.
		\end{equation*}				
		Here  $(\gamma_n)_{n\geq0}$ are the coefficients of  the recurrence relation \eqref{Tre} satisfied by $(S_n(x))_{n\geq 0}.$
	\end{pro}
	Note that the above implies that the coefficients of the  recurrence relation satisfy\begin{equation}\label{relationsgamma}
		\begin{aligned}
			&b_n=(\gamma_{2n+1}+\gamma_{2n}), \ \ &n\geq0, \ \ \  &a_{n}=\gamma_{2n}\gamma_{2n-1},\ \ &n\geq1,\\
			&d_n=(\gamma_{2n+2}+\gamma_{2n+1}),  \ \ &n\geq0, \ \ \ &c_{n}=\gamma_{2n+1}\gamma_{2n},\ \ &n\geq1,
		\end{aligned}
	\end{equation}
	with the convention  $\gamma_0=0$ (see \cite[Chapter 1, Section 9]{GMM21}).
	
	\section{Truncated Laguerre linear functional and the symmetrization process}\label{section3}
	Let $(P_n)_{n\geq0}$ be the sequence of monic orthogonal polynomials with respect to the linear functional
	\begin{equation}\label{laguerrefunc}
		\langle\elln,p(x)\rangle=\int_0^{z} p(x)\,x^{\alpha} e^{-x}dx, \quad p(x)\in\mathbb{P},
	\end{equation}
	with $\alpha>-1$ and $z>0$.
	This functional is known in the literature as \textit{the truncated Gamma functional} and the polynomials $ (P_n)_{n\geq 0}$ are called truncated Laguerre orthogonal polynomials.
	Note that the moments of the functional $\elln$ are given by
	$$\elln_m=\int_0^z x^{m+\alpha}e^{-x}dx=\widehat\gamma(m+\alpha+1,z)$$
	where $\widehat\gamma(a,z)$ is the incomplete gamma function defined by \cite[8.2.1]{OLB10}
	$$\widehat\gamma(a,z)=\int_0^z x^{a-1}e^{-x}dx. $$ Its  series representation , see \cite[8.5.1]{OLB10}, is
	$$\widehat\gamma(a,z)=\dfrac{z^ae^{-z}}{a}\sum_{k=0}^\infty\dfrac{z^k}{k(a+1)_k}$$
	where $(a)_k$ is the Pochhammer symbol defined by $(a)_0=1,$ $(a)_k=a(a+1)\ldots (a+k-1),$ for $k=1,2,\ldots$. In particular,
	$$\elln_m=\dfrac{z^{m+\alpha+1}e^{-z}}{(m+\alpha+1)}\sum_{k=0}^\infty\dfrac{z^k}{k(m+\alpha+1)_k}.$$
	and
	$$\elln_{m+1}=(m+\alpha+1)\elln_m-z^{m+\alpha+1} e^{-z},\quad m\geq 0.$$
	\begin{remark}
		Taking into account the asymptotic expansion  \cite[8.11.2]{OLB10}
		$$\widehat\gamma(a,z)\sim \Gamma(a)\left[1-z^{a-1}e^{-z}\sum_{k=0}^\infty\dfrac{z^{-k}}{\Gamma(a-k)}\right],$$
		then
		\begin{equation}\label{asymp}
			\dfrac{\elln_{m}}{\Gamma(\alpha+m+1)}\sim 1-e^{-z}\sum_{k=0}^\infty\dfrac{z^{\alpha+m-k}}{\Gamma(a-k)} \quad z\to \infty.
		\end{equation}
		Notice that $\Gamma(\alpha+m+1)$ is just the $m$-th moment associated with the Laguerre weight of parameter~$\alpha$.
	\end{remark}
	\begin{teo}
		The functional $\elln$ is semi-classical of class one. Moreover, $\elln$ satisfies the Pearson equation \begin{equation}\label{pearson_eq}
			D(\phi\elln)+\psi\elln=0,
		\end{equation}
		with $\phi(x)=(x-z)x$ and $\psi(x)=-x(z+2+\alpha-x)+z(1+\alpha).$
	\end{teo}
	\begin{proof} Using integration by parts
		\begin{align*}
			\prodint{D(\phi(x)\elln),p(x)}&=-\prodint{\elln,\phi(x) p^\prime(x)}=-\int_0^z\phi(x) p^\prime(x)x^\alpha e^{-x}dx\\
			&=-p(x)\phi(x)x^\alpha e^{-x}\left.\right|_{0}^{z}+\int_0^z p(x)\left(\phi^\prime(x)+\alpha\dfrac{\phi(x)}{x}-\phi(x)\right)x^\alpha e^{-x}dx \\
			&=-\prodint{\elln,\psi(x)p(x)}.
		\end{align*}
	\end{proof}
	\begin{pro}
		Let $(\elln_n)_{n\geq 0}$ the sequence of moments of $\elln$. They satisfy the recurrence relation
		\begin{equation}\label{recmoments}
			\elln_{n+2}-(n+z+2+\alpha)\elln_{n+1}+z(n+\alpha+1)\elln_n=0, \quad n\geq0,
		\end{equation}
		with initial conditions
		$$\elln_0=\widehat\gamma(\alpha+1,z),\quad \elln_1=(\alpha+1)\elln_0-z^{\alpha+1}e^{-z}. $$
	\end{pro}
	\begin{proof}
		From the Pearson equation \eqref{pearson_eq} we have
		\begin{align*}
			\prodint{D(\phi\elln)+\psi\elln,x^n}&=-\prodint{\elln,n\phi(x)\,x^{n-1}}+\prodint{\elln,\psi(x)\,x^n}\\
			&=-\prodint{\elln,n(x-z)x^n}+\prodint{\elln,-x^{n+1}(z+\alpha+2-x)+z(1+\alpha)x^n}\\
			&=-n\elln_{n+1}+nz\elln_{n}-(z+\alpha+2)\elln_{n+1}+\elln_{n+2}+z(1+\alpha)\elln_n\\
			&=\elln_{n+2}-(n+z+\alpha+2)\elln_{n+1}+z(n+\alpha+1)\elln_n=0.
		\end{align*}
	\end{proof}
	\begin{pro}\label{Stielt}
		The Stieltjes function $\Su_{\elln}(t;z)$ associated with the linear functional $\elln$ satisfies the first order non-homogeneous ODE
		$$\phi(t)\partial_t\Su_{\elln}(t;z)+[\phi^\prime(t)+\psi(t)]\Su_{\elln}(t;z)=(t-z-\alpha-1)\elln_0+\elln_1.$$
	\end{pro}
	\begin{proof}
		First, note that
		$$\partial_t\Su_{\elln}(t;z)=-\sum_{n=0}^\infty(n+1)\dfrac{\elln_n}{t^{n+2}}.$$
		Now, from \eqref{recmoments}
		\begin{align*}
			\Su_{\elln}(t;z)&=\dfrac{\elln_0}{t}+\dfrac{\elln_1}{t^2}+\sum_{n=2}^\infty\dfrac{\elln_n}{t^{n+1}}\\
			&=\dfrac{\elln_0}{t}+\dfrac{\elln_1}{t^2}+\sum_{n=1}^\infty(n+1+z+\alpha)\dfrac{\elln_n}{t^{n+2}}-z\sum_{n=1}^\infty(n+1+\alpha)\dfrac{\elln_n}{t^{n+3}} \\
			&=\dfrac{\elln_0}{t}+\dfrac{\elln_1}{t^2}+\dfrac{(z+\alpha)}{t}\left(\sum_{n=1}^\infty\dfrac{\elln_n}{t^{n+1}}+\dfrac{\elln_0}{t}-\dfrac{\elln_0}{t}\right)+\left(\sum_{n=1}^\infty(n+1)\dfrac{\elln_n}{t^{n+2}}+\dfrac{\elln_0}{t^2}-\dfrac{\elln_0}{t^2}\right)\\
			&-\dfrac{z\alpha}{t^2}\sum_{n=0}^\infty\dfrac{\elln_n}{t^{n+1}}-\dfrac{z}{t}\sum_{n=0}^\infty(n+1)\dfrac{\elln_n}{t^{n+1}}\\
			&=\dfrac{\elln_0}{t}+\dfrac{\elln_1}{t^2}+\dfrac{(z+\alpha)}{t}\left(\Su_{\elln}(t;z)-\dfrac{\elln_0}{t}\right)-\partial_t\Su_{\elln}(t;z)-\dfrac{\elln_0}{t^2}-\dfrac{z\alpha}{t^2}\Su_{\elln}(t;z)+\dfrac{z}{t}\partial_t\Su_{\elln}(t;z).
		\end{align*}
		Thus, multiplying both sides of the equation by $t^2$ and reorganizing
		$$(t^2-zt)\partial_t\Su_{\elln}(t;z)+(t^2-(z+\alpha)t+z\alpha)\Su_{\elln}(t;z)=(t-z-\alpha-1)\elln_0+\elln_1.$$
	\end{proof}
	\begin{teo}\label{sympol}
		Let $\widetilde\elln$ be the linear functional defined by $\widetilde\elln=x\elln$, then \begin{enumerate}
			\item $\widetilde\elln$ is a quasi-definite functional as well as  semi-classical of class one.
			\item If $(Q_n)_{n\geq 0}$ is the sequence of monic orthogonal polynomials with respect to $\widetilde\elln$, then
			$$xQ_n(x)=P_{n+1}(x)-\dfrac{P_{n+1}(0)}{P_n(0)}P_n(x), \quad n\geq 0.$$
		\end{enumerate}
	\end{teo}
	Let $(P_n)_{n\geq 0}$ and $(\widetilde P_n)_{n\geq 0}$ be the MOPS with respect to $\elln$ and $\widetilde\elln$, respectively. Define
	\begin{equation}\label{simetricpoly}
		S_{2n}(x)=P_n(x^2) \ \ \ \ \ \text{and}\ \ \ \ \ S_{2n+1}(x)=xQ_n(x^2).
	\end{equation}
	As a consequence of the above,
	\begin{pro} Let $\elln$ be the linear functional defined in \eqref{laguerrefunc} with parameter  $\alpha> -1$  and    $(S_n)_{n\geq0}$  be the polynomials defined in \eqref{simetricpoly}, Then
		\begin{enumerate}
			\item $(S_n)_{n\geq 0}$ is the sequence of orthogonal polynomials with respect to the symmetric linear functional
			\begin{equation*}
				\prodint{\hn,p}= \int_{-\sqrt{z}}^{\sqrt{z}} p(x)\,|x|^{2\alpha+1} e^{-x^2}dx,  \quad \alpha >-1, \quad p(x)\in\mathbb{P}.
			\end{equation*}
			\item The moments $(\hn_n)_{n\geq0}$ associated with the functional $\hn$ satisfy
			$$\hn_{2n+4}-(n+z+2+\alpha)\hn_{2n+2}+z(n+\alpha+1)\hn_{2n}=0,\quad n\geq 0,$$
			with initial conditions
			$$\hn_0=\widehat\gamma(\alpha+1,z),\quad \hn_2=(\alpha+1)\hn_0-z^{\alpha+1}e^{-z}. $$
			\item The following  ratio asymptotic holds $$\dfrac{\hn_{2m}}{\Gamma(\alpha+m+1)}\sim 1-e^{-z}\sum_{k=0}^\infty\dfrac{z^{\alpha+m-k}}{\Gamma(a-k)},\quad z\to \infty.$$
			\item The functional $\hn$ is semiclassical of class 2 if $\alpha= -1/2$ and of class 3 in the other case. Moreover, $\hn$ satisfies the Pearson equation \begin{equation*}
				D(\Phi\hn)+\Psi\hn=0,
			\end{equation*}
			with
			\begin{equation*}
				\begin{cases}
					\Phi(x)=(x^2-z)x \quad  \Psi(x)=2x^4-2(2+\alpha+z)x^2+2z(\alpha+1), &  \alpha\ne -1/2, \\[10pt]
					\Phi(x)=(x^2-z), \quad  \Psi(x)=2x(x^2-(z+1)),&  \alpha= -1/2.
				\end{cases}
			\end{equation*}
			\item The Stieltjes function $\Su_{\hn}(t;z)$ associated with the linear functional $\hn$ satisfies the first-order non-homogeneous ODE
			\begin{multline*}
				t(t^2-z)\,\partial_t\Su_{\hn}(t;z)+2\left[t^4-(z+\alpha+1/2)t^2+z(\alpha+1/2)\right]\Su_{\hn}(t;z)=2t[t^2-(z+\alpha+1)]\hn_0+2t\hn_2.
			\end{multline*}
		\end{enumerate}
	\end{pro}
	\begin{proof}
		Since $|x|^{2\alpha+1} e^{-x^2}$ is an even function, it is clear that $\prodint{\hn,S_{2n}S_{2m+1}}=0,$ where the polynomials $(S_n)_{n\geq0}$ were defined in \eqref{simetricpoly}. Moreover,
		\begin{align*}
			\prodint{\hn,x^{2m} p(x^2)q(x^2)}=\int_{0}^{z} p(x)q(x)x^{\alpha+m}e^{-x}dx,\quad m\geq 0,,
		\end{align*}
		and the orthogonality is easily obtained from Theorem \ref{sympol}. Properties (ii) and (iii)  are a direct consequence of \eqref{recmoments} and \eqref{asymp}, respectively. We will show (iv) for $\alpha\ne-1/2$ the other case is similar. The Pearson equation is a direct consequence of integration by parts. To check the class of functional, we will use Proposition \ref{sim_cond}.  Since in our case \eqref{prodformula} becomes
		$$\prod_{c:\,\phi(c)=0}\left(\left|(c^2 - z)(-c^2 + \alpha + 1/2)\right|+\left|\,c\,[c^2-(\alpha+z+1)]\hn_0+c\hn_2\right|\right),$$
		with $c=0$, $c=\pm\sqrt{z},$ the zeros of $\phi$, then the product will be different to zero. Therefore, the class of $\hn$ is $3$. Finally, (v) can be deduced as in Proposition \ref{Stielt}.
	\end{proof}
	\begin{remark}\label{rem3,7}
		The case where $\alpha=-1/2$ was studied by D. Dominici and F. Marcellán in \cite{DM22}. Thus, hereinafter we shall deal with the case $\alpha\ne-1/2.$ On the other hand, when $\mu= \alpha+1/2 $  and $z= \infty$ you get the so called generalized Hermite polynomials analyzed in \cite{Ted} (see also \cite[pp.156-158]{Chi}, \cite{Gabor}). These generalized Hermite  polynomials as well as the Laguerre polynomials are Brenke orthogonal polynomials. In fact, the corresponding generating function $\sum_{n=0}^{\infty} p_{n}(x) w^{n}$ has the form $A(w) B(xw),$ where $A(w)= \sum_{n=0}^{\infty} a_{n} w^{n},$ $B(w)= \sum_{n=0}^{\infty} b_{n} w^{n},$ In such a way,  $p_{n}(x)= \sum_{k=0}^{n} a_{n-k} b_{k}x^{k},$ The sequences of orthogonal polynomials which are also Brenke polynomials were determined in \cite{Chihara}.
	\end{remark}
	
	\section{Laguerre-Freud's equations}\label{section4}
	If $\elln$ is a semiclassical linear functional satisfying \eqref{pearson-semic}, then the coefficients of the three-term recurrence relation \eqref{ttrr} satisfy a nonlinear system of equations obtained from the relations.
	\begin{align}
		\label{1eq}&\prodint{\psi\elln,P_n^2}=-\prodint{D(\phi\elln),P_n^2},\\
		\label{2eq}&\prodint{\psi\elln,P_{n+1}P_n}=-\prodint{D(\phi\elln),P_{n+1}P_n}.
	\end{align}
	They are known in the literature as \textit{Laguerre-Freud equations} (see \cite{Said}). Let $(P_n)_{n\geq 0}$ be the sequence of monic orthogonal polynomials with respect to $\elln$ that satisfies the recurrence relation \eqref{ttrr}. If
	\begin{equation}\label{rev}
		P_n(x)=x^n+\sum_{k=0}^{n-1}\lambda_{n,k}x^k,
	\end{equation}
	then
	\begin{align*}
		xP_n(x)&=x^{n+1}+\lambda_{n,n-1}x^{n}+\lambda_{n,n-2}x^{n-1}+\lambda_{n,n-3}x^{n-2}+\mathcal{O}(x^{n-3})\\
		&=x^{n+1}+(\lambda_{n+1,n}+b_n)x^n+(\lambda_{n+1,n-1}+\lambda_{n,n-1}b_n+a_n)x^{n-1}+\mathcal{O}(x^{n-2}).
	\end{align*}
	In particular
	\begin{equation}\label{cases}
		\begin{cases}
			b_n=\lambda_{n,n-1}-\lambda_{n+1,n},  \\
			\lambda_{n,n-1}b_n+a_n=\lambda_{n,n-2}-\lambda_{n+1,n-1}.
		\end{cases}
	\end{equation}
	Now, reverse \eqref{rev} and take into account \eqref{cases}
	\begin{align}
		\notag\partial_xP_n(x)&=nx^{n-1}+(n-1)\lambda_{n,n-1}x^{n-2}+(n-2)\lambda_{n,n-2}x^{n-3}+\mathcal{O}(x^{n-4})\\
		\notag&=nP_{n-1}(x)+((n-1)\lambda_{n,n-1}-n\lambda_{n-1,n-2})P_{n-2}(x)\\
		\label{der}&+(n\lambda_{n-1,n-2}\lambda_{n-2,n-3}-n\lambda_{n-1,n-3}-(n-1)\lambda_{n,n-1}\lambda_{n-2,n-3}\\
		\notag&+(n-2)\lambda_{n,n-2})P_{n-3}+\mathcal{O}(x^{n-4})
		\\ \notag&=nP_{n-1}(x)-(nb_{n-1}+\lambda_{n,n-1})P_{n-2}(x)+\\
		\notag &+(n(b_{n-1}b_{n-2}-a_{n-1})-2\lambda_{n,n-2}+\lambda_{n,n-1}\lambda_{n-2,n-3})P_{n-3}+\mathcal{O}(x^{n-4}).
	\end{align}
	Keeping this in mind,
	\begin{pro}
		The coefficients of the three term  recurrence relation \eqref{ttrr} associated with the linear functional $\elln$ satisfy the Laguerre-Freud equations
		\begin{equation}\label{111}
			a_{n+1}-a_{n-1}+b_n^2-b_{n-1}^2+(z+\alpha+2n)(b_{n-1}-b_{n})-2(b_{n-1}+b_{n})+2z=0.
		\end{equation}
		\begin{equation}\label{112}
			(2n-b_n+z+\alpha+2)(a_{n+1}-a_{n})+(1-b_{n+1})a_{n+1}+(3+b_{n-1})a_n+b^2_n-zb_n=0.
		\end{equation}
		Moreover, if we define
		\begin{equation}\label{wn}
			\omega_{n}=b_n-\dfrac{\alpha}{2}-\dfrac{z}{2}-n,
		\end{equation}
		the expressions \eqref{111} and \eqref{112} become
		\begin{align}
			a_{n+1}-a_{n-1}+(\omega_{n}+\omega_{n-1}-1)(\omega_{n}-\omega_{n-1}-1)=2(\alpha+2n),\label{w1}\\
			a_n(\omega_{n}+\omega_{n-1})-a_{n+1}(\omega_{n}+\omega_{n+1}-2)+\left(\omega_n+\dfrac{\alpha}{2}+n\right)^2=\dfrac{z^2}{4}.\label{w2}
		\end{align}
	\end{pro}
	\begin{proof}
		Let $h_n=\prodint{\elln,P_n^2}$.
		Equation \eqref{1eq} is equivalent to
		\begin{equation}\label{eqq1}
			\prodint{\elln,\psi P_n^2}=2\prodint{\elln,\phi P_nP^\prime_n},
		\end{equation}
		where $\phi(x)=(x-z)x$ and $\psi(x)=-x(z+2+\alpha-x)+z(1+\alpha).$   Define $d=z+2+\alpha$ and $e=z(1+\alpha)$. The left-hand side of \eqref{eqq1}
		\begin{align*}
			\prodint{\elln,(x^2-dx+e)P_nP_n} &=\prodint{\elln,x(P_{n+1}+b_nP_n+a_{n}P_{n-1})P_n}\\&-d\prodint{\elln,(P_{n+1}+b_nP_n+a_{n}P_{n-1})P_n}+eh_n\\
			&=(a_{n+1}+b_n^2+a_n-db_n+e)h_n.
		\end{align*}
		On the other hand, using \eqref{der}, the right-hand side of \eqref{eqq1} reads
		\begin{equation}\label{ppp}
			\begin{aligned}
				&\prodint{\elln ,\phi P_nP^\prime_n}=\prodint{\elln, P_n\left[n\phi P_{n-1}-(nb_{n-1}+\lambda_{n,n-1})\phi P_{n-2}+\mathcal{O}(x^{n-1})\right]}\\
				&=n\prodint{\elln,xP_{n}(P_n+b_{n-1}P_{n-1}+a_{n-1}P_{n-2})}-z\,n\prodint{\elln,P_n(P_n+b_{n-1}P_{n-1}+a_{n-1}P_{n-2})}\\
				&-(nb_{n-1}+\lambda_{n,n-1})\prodint{\elln,xP_nP_{n-1}}\\
				&=\left[nb_n+nb_{n-1}-zn-(nb_{n-1}+\lambda_{n,n-1})\right]h_n.
			\end{aligned}
		\end{equation}
		From the above,
		\begin{equation}\label{11}
			a_{n+1}+b_{n}^2+a_{n}-db_{n}+e=2(nb_n-nz-\lambda_{n,n-1})
		\end{equation}
		shifting $n\to n-1$
		\begin{equation}\label{12}
			a_{n}+b_{n-1}^2+a_{n-1}-db_{n-1}+e=2((n-1)b_{n-1}-(n-1)z-\lambda_{n-1,n-2}).
		\end{equation}
		Subtracting \eqref{11} from \eqref{12} and using \eqref{cases}, we get \eqref{111}. On the other hand, \eqref{2eq} is equivalent to
		\begin{equation*}
			\prodint{\elln,\psi P_nP_{n+1}}=\prodint{\elln,\phi P_nP^\prime_{n+1}}+\prodint{\elln,\phi P^\prime_{n}P_{n+1}}.   \end{equation*}
		Now
		\begin{equation}\label{eqq2}
			\begin{aligned}
				\prodint{\elln,\psi P_nP_{n+1}}&=\left[b_{n+1}+b_n-d\right]a_{n+1}h_n,\\
				\prodint{\elln,\phi P^\prime_{n}P_{n+1}}&=na_{n+1}h_n,
				\\\prodint{\elln,\phi P_nP^\prime_{n+1}}&=\left[(n+1)a_{n+1}-\lambda_{n+1,n}(b_n+b_{n-1}-\lambda_{n-1,n-2})+z\lambda_{n+1,n}-2\lambda_{n+1,n-1}\right]h_n.
			\end{aligned}
		\end{equation}
		From the above and \eqref{cases}, we have the following.
		\begin{equation*}
			(2n+1)a_{n+1}-\lambda_{n+1,n}(b_n+b_{n-1}-\lambda_{n-1,n-2})+z\lambda_{n+1,n}-2\lambda_{n+1,n-1}=  \left[b_{n+1}+b_n-d\right]a_{n+1}
		\end{equation*}
		
		shifting $n\to n-1$, subtracting and taking into account \eqref{cases} \begin{align*}
			2n(a_{n+1}-a_n)&+a_{n+1}+a_n+b_{n}b_{n-1}-\lambda_{n+1,n}(b_n-\lambda_{n-1,n-2})+\lambda_{n,n-1}(b_{n-2}-\lambda_{n-2,n-3}) \\&-zb_n+2(\lambda_{n,n-1}b_n+a_n)\\
			&=(b_n-d)(a_{n+1}-a_n)+b_{n+1}a_{n+1}-b_{n-1}a_n.
		\end{align*}
		On the other hand,
		\begin{equation}\label{lafreu2}
			\begin{aligned}
				-&\lambda_{n+1,n}(b_n-\lambda_{n-1,n-2})+\lambda_{n,n-1}(b_{n-2}-\lambda_{n-2,n-3})+2\lambda_{n,n-1}b_n\\
				=&-(\lambda_{n,n-1}-b_n)(b_n-\lambda_{n-1,n-2})-\lambda_{n,n-1}\lambda_{n-1,n-2}+2\lambda_{n,n-1}b_n\\
				=&b_n^2-b_{n}b_{n-1}.
			\end{aligned}
		\end{equation}
		Thus, replacing the above in \eqref{lafreu2}  we obtain \eqref{112}.
		
		Finally, taking into account that
		\begin{multline*}b_n^2-b_{n-1}^2+(z+\alpha+2n)(b_{n-1}-b_{n})-2(b_{n-1}+b_{n})\\=(b_n+b_{n-1}-2)(b_n+b_{n-1}-[z+\alpha+2n])-2(z+\alpha+2n)\end{multline*} substituting this in \eqref{111} and \eqref{112},  and defining $w_n$ as in \eqref{wn} we obtain \eqref{w1} and \eqref{w2}, respectively.
	\end{proof}
	
	\begin{remark}
		Note that the coefficients of the recurrence relation \eqref{ttrr} associated with the linear functional $\widetilde\elln$ satisfy the Laguerre-Freud equations
		\begin{equation}\label{222}
			c_{n+1}-c_{n-1}+d_n^2-d_{n-1}^2+(z+\alpha+1+2n)(d_{n-1}-d_{n})-2(d_{n-1}+d_{n})+2z=0.
		\end{equation}
		\begin{equation}\label{223}
			(2n-d_n+z+\alpha+3)(c_{n+1}-c_{n})+(1-d_{n+1})c_{n+1}+(3+d_{n-1})c_n+d^2_n-zd_n=0.
		\end{equation}
		Defining
		\begin{equation*}
			\widetilde{\omega}_{n}=d_n-\dfrac{\alpha+1}{2}-\dfrac{z}{2}-n,   \end{equation*}
		the expressions \eqref{222} and \eqref{223} become
		\begin{align*}
			c_{n+1}-c_{n-1}+(\widetilde{\omega}_{n}+\widetilde{\omega}_{n-1}-1)(\widetilde{\omega}_{n}-\widetilde{\omega}_{n-1}-1)=2(\alpha+2n+1),\\
			c_n(\widetilde{\omega}_{n}+\widetilde{\omega}_{n-1})-c_{n+1}(\widetilde{\omega}_{n}+\widetilde{\omega}_{n+1}-2)+\left(\widetilde{\omega}_n+\dfrac{\alpha+1}{2}+n\right)^2=\dfrac{z^2}{4}.
		\end{align*}
		
	\end{remark}
	
	\begin{coro}
		The coefficients of the three term recurrence relation \eqref{Tre} associated with the linear functional $\hn$ satisfy
		\begin{multline*}
			\gamma_{n+2}\,\gamma_{n+1}+(\gamma_{n+1}+\gamma_{n})^2-\left(z+\alpha+n+2\right)(\gamma_{n+1}+\gamma_{n})\\
			-\gamma_{n-2}\,\gamma_{n-3}-(\gamma_{n-1}+\gamma_{n-2})^2+\left(z+\alpha+n-2\right)(\gamma_{n-1}+\gamma_{n-2})=-2z.
		\end{multline*}
		\begin{multline*}
			\left(n+z+\alpha+3-\sum_{k=n}^{n+3}\gamma_k\right)\gamma_{n+2}\gamma_{n+1}+\left(n+z+\alpha+1-\sum_{k=n-2}^{n+1}\gamma_k\right)\gamma_{n}\gamma_{n-1}+\\+(\gamma_{n+1}+\gamma_{n})(\gamma_{n+1}+\gamma_{n}-z)=0.
		\end{multline*}
	\end{coro}
	\begin{proof} The result is a direct consequence of \eqref{111}, \eqref{112}, \eqref{222}, \eqref{223} and \eqref{relationsgamma}.
	\end{proof}
	
	\begin{teo}\label{Teo4.4} Let $(a_n)_{n\geq 1}$, $(b_n)_{n\geq 0}$ be the coefficients of the recurrence relation \eqref{ttrr}. Then,   they satisfy the following nonlinear  relation
 \begin{multline*}
 (\alpha-b_{n}+2n+1)(\alpha-b_{n-1}+2n-1)a_n\\=\dfrac{1}{4}\left[a_n-a_{n+1}-b_n\left(b_n-[z+\alpha+2(n+1)]\right)-(1+\alpha+2n)z\right]^2.
 \end{multline*}
 	\end{teo}
	\begin{proof}
		Using integration by parts, equation \eqref{der}, and the fact that $a_{n-1}=h_{n-1}/h_{n-2}$ we get
		$$
		\begin{aligned}
			a_nh_{n-1}&=\prodint{\elln,xP_nP_{n-1}}=-P_n(z)P_{n-1}(z)z^{\alpha+1}e^{-z}+\prodint{\elln,xP^\prime_{n}P_{n-1}}\\&=-P_n(z)P_{n-1}(z)z^{\alpha+1}e^{-z}+\prodint{\elln,(P_n+b_{n-1}P_{n-1}+a_{n-1}P_{n-2})P_n^\prime}\\
			&=-P_n(z)P_{n-1}(z)z^{\alpha+1}e^{-z}+nb_{n-1}h_{n-1}-a_{n-1}\prodint{\elln,P_{n-2}P_n^\prime}\\
			&=-P_n(z)P_{n-1}(z)z^{\alpha+1}e^{-z}-h_{n-1}\lambda_{n,n-1}.
		\end{aligned}
		$$
		Thus, from the above and \eqref{11}
		\begin{equation}\label{nonlinearrecubnan}
			\begin{aligned}
				P^2_n(z)P^2_{n-1}(z)z^{2\alpha+2}e^{-2z}&=(a_n+\lambda_{n,n-1})^2h_{n-1}^2\\
				&=\dfrac{1}{4}\left[a_n-a_{n+1}-b_n\left(b_n-[z+\alpha+2(n+1)]\right)-(1+\alpha+2n)z\right]^2.
			\end{aligned}
		\end{equation}
		On the other hand
		$$\begin{aligned}
			b_nh_n=\prodint{\elln,xP_n^2}&=-P^2_n(z)z^{\alpha+1}e^{-z}+(\alpha+1)h_n+2\prodint{\elln,xP_nP_n^\prime}\\&=-P^2_n(z)z^{\alpha+1}e^{-z}+(\alpha+1)h_n+2nh_n.
		\end{aligned}$$
		Thus,
		$$P^2_n(z)z^{\alpha+1}e^{-z}=(\alpha+1+2n-b_n)h_n.$$
		Replacing  this formula in \eqref{nonlinearrecubnan} and taking into account that $a_n=h_n/h_{n-1}$, we get the result.
	\end{proof}
	\begin{pro}
		The coefficients of the three term  recurrence relation \eqref{Tre} associated with the linear functional $\hn$ satisfy the Laguerre-Freud equation
		\begin{multline}\label{LFh}
			\gamma_{n+1}(\gamma_{n+2}+\gamma_{n+1}+\gamma_{n}-(\alpha+z+n+2))-\gamma_{n-1}(\gamma_{n-2}+\gamma_{n-1}+\gamma_{n}-(\alpha+z+n-1))-\gamma_n+z=0
		\end{multline}
	\end{pro}
	\begin{proof}
		Let $\mathfrak{h}_n=\prodint{\hn,S_n^2}$. Since $\hn$  satisfies \eqref{pearson-semic} with $ \Phi(x)=(x^2-z)x$ and $  \Psi(x)=2x^4-2(2+\alpha+z)x^2+2z(\alpha+1)$, then
		\begin{equation}\label{eqq1h}
			\dfrac{1}{2}\prodint{\hn,\Psi S_n^2}=\prodint{\hn,\Phi S_nS^\prime_n},
		\end{equation}
		Define $d=(2+\alpha+z)$ and $c=z(\alpha+1)$. On the one hand, using \eqref{Tre} and the fact that
		\begin{equation}\label{x2S}
			x^2S_n(x)=S_{n+2}(x)+(\gamma_{n+1}+\gamma_{n})S_{n}(x)+(\gamma_{n-1}\gamma_{n})S_{n-2}(x),
		\end{equation}
		we have
		\begin{equation*}
			\begin{aligned}
				\dfrac{1}{2}&\prodint{\hn,\Psi S_n^2}=\prodint{\hn,(x^4-dx^2+c)S_n^2(x)} \\
				&=\prodint{\hn,x^2S_n(S_{n+2}+(\gamma_{n+1}+\gamma_{n})S_{n}+\gamma_{n-1}\gamma_{n}S_{n-2})}
				-d\prodint{\hn,(\gamma_{n+1}+\gamma_{n})S_n^2}+c\,\mathfrak{h}_n \\
				&=\prodint{\hn,\gamma_{n+1}\gamma_{n+2}S^2_{n}+(\gamma_{n+1}+\gamma_{n})^2S^2_{n}+\gamma_{n-1}\gamma_{n}S^2_{n}}
				-d\prodint{\hn,(\gamma_{n+1}+\gamma_{n})S_n^2}+c\,\mathfrak{h}_n\\&
				=\left(\gamma_{n+1}\gamma_{n+2}+(\gamma_{n+1}+\gamma_{n})^2+\gamma_{n-1}\gamma_{n}-d(\gamma_{n+1}+\gamma_{n})+c\right)\mathfrak{h}_n.
			\end{aligned}
		\end{equation*}
		On the other hand, if
		$S_n(x)=x^n+\widehat{\lambda}_{n,n-2}x^{n-2}+\widehat{\lambda}_{n,n-3}x^{n-3}+\cdots,$
		then from the above, \eqref{Tre} and \eqref{x2S} we get
		$$
		\begin{aligned}
			xS_{n}&=x^{n+1}+(\widehat{\lambda}_{n+1,n-1}+\gamma_n)x^{n-1}+\mathcal{O}(x^{n-3}),\\
			x^2S_n(x)&=x^{n+2}+[\widehat{\lambda}_{n+2,n}+(\gamma_{n+1}+\gamma_{n})]x^n+\mathcal{O}(x^{n-2})
		\end{aligned}
		$$
		In particular,
		\begin{equation*}
			\begin{cases}
				\gamma_{n-1}=\widehat{\lambda}_{n-1,n-3}-\widehat{\lambda}_{n,n-2},  \\
				\gamma_n+\gamma_{n+1}=\widehat{\lambda}_{n,n-2}-\widehat{\lambda}_{n+2,n}.
			\end{cases}
		\end{equation*}
		Moreover,  since
		$$
		\begin{aligned}
			\partial_xS_n&=nS_{n-1}+((n-2)\widehat{\lambda}_{n,n-2}-n\widehat{\lambda}_{n-1,n-3})x^{n-3}+\mathcal{O}(x^{n-5})\\&=nS_{n-1}-(n\gamma_{n-1}+2\widehat{\lambda}_{n,n-2})x^{n-3}+\mathcal{O}(x^{n-5}).
		\end{aligned}
		$$
		This yields
		\begin{equation*}
			\begin{aligned}
				&\prodint{\hn,\Phi S_nS^\prime_n}=\prodint{\hn,x^3 S_nS^\prime_n}-z\prodint{\hn,x S_nS^\prime_n} \\
				&=n\prodint{\hn,xS_n[S_{n+1}+(\gamma_{n}+\gamma_{n-1})S_{n-1}+\gamma_{n-1}\gamma_{n-2}S_{n-2}]}-(n\gamma_{n-1}+2\widehat{\lambda}_{n,n-2}+nz)\mathfrak{h}_n\\
				&=\left[n(\gamma_{n+1}+\gamma_{n})  -2\widehat{\lambda}_{n,n-2}-zn\right]\mathfrak{h}_n.
			\end{aligned}
		\end{equation*}
		Thus, from \eqref{eqq1h} we obtain
		\begin{equation}\label{lambdahat}\gamma_{n+1}\gamma_{n+2}+(\gamma_{n+1}+\gamma_{n})^2+\gamma_{n-1}\gamma_{n}
			-(d+n)(\gamma_{n+1}+\gamma_{n})+c=  -2\widehat{\lambda}_{n,n-2}-nz.\end{equation}
		Shifting $n\to n-1$ and subtracting, we get \eqref{LFh}.
		
	\end{proof}
	
	\begin{teo} Let $(\gamma_n)_{n\geq 1}$ be the coefficients of the three term recurrence relation \eqref{Tre}. They satisfy the following nonlinear  relation
		\begin{equation}\label{nonlinearrec_gamma1}
			\Big((\alpha+n+1)-(\gamma_{n+1}+\gamma_{n})\Big)\Big((\alpha+n-1)-(\gamma_{n-1}+\gamma_{n-2})\Big)\gamma_{n}\gamma_{n-1}=\dfrac{1}{4}(2\gamma_{n}\gamma_{n-1}+2\widehat{\lambda}_{n,n-2})^2.
		\end{equation}
		Moreover,
		\begin{enumerate}
			\item If $n=2k$, we get
			\begin{equation}\label{equality1}
				{\gamma_{2k}}\Big((\alpha+2k+1)-(\gamma_{2k+1}+\gamma_{2k})\Big)\Big((\alpha+2k)-(\gamma_{2k}+\gamma_{2k-1})\Big)=z(k-\gamma_{2k})^2.
			\end{equation}
			\item If $n=2k+1$, we get
			\begin{equation}\label{equality2}
				\gamma_{2k+1}\Big((\alpha+2k+2)-(\gamma_{2k+2}+\gamma_{2k+1})\Big)\Big((\alpha+2k+1)-(\gamma_{2k+1}+\gamma_{2k})\Big)=z(k+\alpha+1-\gamma_{2k+1})^2.
			\end{equation}
		\end{enumerate}
	\end{teo}
	\begin{proof}
		Notice that
		$$\prodint{\hn,\partial_{x}(xS^2_n(x))}=\prodint{\hn,2xS_nS_n^\prime}+\prodint{\hn,S_n^2}=(2n+1)\mathfrak{h}_n.$$
		On the other hand
		$$
		\begin{aligned}
			\prodint{\hn,\partial_{x}(xS^2_n(x))}&=2z^{\alpha+1}S^2_{n}\left(\sqrt{z}\right)e^{-z}-(2\alpha+1)\prodint{\hn,S_n^2}+2\prodint{\hn,x^2S_n^2}\\
			&=2z^{\alpha+1}S^2_{n}\left(\sqrt{z}\right)e^{-z}-(2\alpha+1)\mathfrak{h}_n+2(\gamma_{n+1}+\gamma_n)\mathfrak{h}_n.
		\end{aligned}
		$$
		From the above we obtain
		\begin{equation}\label{exhn}
			z^{\alpha+1}S^2_{n}\left(\sqrt{z}\right)e^{-z}=\left[(\alpha+n+1)-(\gamma_{n+1}+\gamma_n)\right]\mathfrak{h}_n.
		\end{equation}
		Now observe that
		$$\begin{aligned}
			\gamma_{n-1}\gamma_{n}\mathfrak{h}_{n-2}&=\prodint{\hn,x^2S_{n}S_{n-2}}=-z^{\alpha+1}S_{n}\left(\sqrt{z}\right)S_{n-2}e^{-z}\left(\sqrt{z}\right)+\dfrac{1}{2}\prodint{\hn,xS_n^\prime S_{n-2}+xS_n S^\prime_{n-2}}   \\
			&=-S_{n}\left(\sqrt{z}\right)S_{n-2}\left(\sqrt{z}\right)z^{\alpha+1}e^{-z}+\dfrac{1}{2}\left(n\gamma_{n-1}-(n\gamma_{n-1}+2\widehat{\lambda}_{n,n-2})\right)\mathfrak{h}_{n-2}.
		\end{aligned}$$
		Hence,
		$$S^2_{n}\left(\sqrt{z}\right)S^2_{n-2}\left(\sqrt{z}\right)z^{2\alpha+2}e^{-2z}=(\gamma_{n-1}\gamma_{n}+\widehat{\lambda}_{n,n-2})^2\mathfrak{h}_{n-2}.$$
		Thus, from \eqref{exhn} we get
		$$\Big((\alpha+n+1)-(\gamma_{n+1}+\gamma_{n})\Big)\Big((\alpha+n-1)-(\gamma_{n-1}+\gamma_{n-2})\Big)\gamma_{n}\gamma_{n-1}=\dfrac{1}{4}(2\gamma_{n}\gamma_{n-1}+2\widehat{\lambda}_{n,n-2})^2$$ and the first part of Theorem is proved.
		Now, taking into account that
		\begin{equation}\label{pareimparSn}
			\begin{aligned}
				\gamma_{n}&\mathfrak{h}_{n-1}=\prodint{\hn,x S_n S_{n-1}}\\
				&=-S_{n}\left(\sqrt{z}\right)S_{n-1}\left(\sqrt{z}\right)z^{\alpha+1/2}e^{-z}+\dfrac{1}{2}\prodint{\hn,S_nS^\prime_{n-1}+S_{n-1}S^\prime_{n}}+\dfrac{(2\alpha+1)}{2}\prodint{\hn,\dfrac{S_{n}S_{n-1}}{x}}\\
				&=-S_{n}\left(\sqrt{z}\right)S_{n-1}\left(\sqrt{z}\right)z^{\alpha+1/2}e^{-z}+\dfrac{n}{2}\mathfrak{h}_{n-1}+\dfrac{(2\alpha+1)}{2}\prodint{\hn,\dfrac{S_{n}S_{n-1}}{x}}
			\end{aligned}
		\end{equation}
		the last term of the above equation is well defined, since from \eqref{pareimpar} if $S_n$ is an even (odd) function, then $S_{n-1}$ is an odd (even) function.
		\begin{itemize}
			
			\item Suppose that $n=2k$,  then \eqref{pareimparSn} becomes
			$$S_{2k}\left(\sqrt{z}\right)S_{2k-1}\left(\sqrt{z}\right)z^{\alpha+1/2}e^{-z}=[k-\gamma_{2k}]\mathfrak{h}_{2k-1}.$$
			Using \eqref{exhn} we get
			$$\dfrac{\gamma_{2k}}{z}\Big((\alpha+2k+1)-(\gamma_{2k+1}+\gamma_{2k})\Big)\Big((\alpha+2k)-(\gamma_{2k}+\gamma_{2k-1})\Big)=(k-\gamma_{2k})^2.$$
			\item Suppose now  that $n=2k+1$,  then \eqref{pareimparSn} became  $$S_{2k+1}\left(\sqrt{z}\right)S_{2k}\left(\sqrt{z}\right)z^{\alpha+1/2}e^{-z}=[k+\alpha+1-\gamma_{2k+1}]\mathfrak{h}_{2k}.$$
			Using \eqref{exhn} we get
			$$\dfrac{\gamma_{2k+1}}{z}\Big((\alpha+2k+2)-(\gamma_{2k+2}+\gamma_{2k+1})\Big)\Big((\alpha+2k+1)-(\gamma_{2k+1}+\gamma_{2k})\Big)=(k+\alpha+1-\gamma_{2k+1})^2.$$
			
		\end{itemize}
	\end{proof}
	\begin{coro}\label{coro_simetricos_painleve}
		If we define
		$g_{n}=\gamma_{n}-\dfrac{\alpha}{2}-\dfrac{n}{2}-\dfrac{1}{4}$, then $(g_n)_{n\geq1}$ satisfies the following nonlinear recurrence relations.
		\begin{multline}\label{gteorem}
			(g_n+g_{n+1})\left(2g_n+\alpha+n+\dfrac{1}{2}\right)(g_{n-2}+g_{n-1})\left(2g_{n-1}+\alpha+n-\dfrac{1}{2}\right)=\\ \Big((g_n+g_{n+1}-z)\left(2g_n+\alpha+n+\dfrac{1}{2}\right)+(g_{n-2}+g_{n-1}-z)\left(2g_{n-1}+\alpha+n-\dfrac{1}{2}\right)+2z(n+\alpha)\Big)^2
		\end{multline}
		as well as
		\begin{equation}\label{painlevesimetricos}
			\left(\dfrac{\alpha}{2}+\dfrac{n}{2}+\dfrac{1}{4}+g_{n}\right)(g_{n}+g_{n+1})(g_{n}+g_{n-1})=z\left(\dfrac{\alpha}{2}+\dfrac{1}{4}-g_{n}\right)^2.
		\end{equation}
		
	\end{coro}
	\begin{proof}
		Note that from \eqref{lambdahat} and \eqref{LFh}
		$$\begin{aligned}
			-2\gamma_{n}\gamma_{n-1}-2\widehat{\lambda}_{n,n-2}=&\gamma_{n+1}(\gamma_{n+2}+\gamma_{n+1}+\gamma_{n}-(\alpha+z+n+2))\\ &+\gamma_{n}(\gamma_{n+1}+\gamma_{n}-\gamma_{n-1}-(\alpha+z+n+2))+z(n+\alpha+1)
			\\=&\gamma_{n}(\gamma_{n+1}+\gamma_{n}-(\alpha+z+n+1))\\ &+\gamma_{n-1}(\gamma_{n-2}+\gamma_{n-1}-(\alpha+z+n-1))+z(n+\alpha).
		\end{aligned}$$
		Replacing the above in \eqref{nonlinearrec_gamma1} and taking into account the definition of $g_n$, we get \eqref{gteorem}.
		Now, from \eqref{equality1} and \eqref{equality2} we get
		\begin{align*}
			&\left(\dfrac{\alpha}{2}+k+\dfrac{1}{4}+g_{2k}\right)(g_{2k}+g_{2k+1})(g_{2k}+g_{2k-1})=z\left(\dfrac{\alpha}{2}+\dfrac{1}{4}-g_{2k}\right)^2,\\
			&\left(\dfrac{\alpha}{2}+k+\dfrac{3}{4}+g_{2k+1}\right)(g_{2k+2}+g_{2k+1})(g_{2k}+g_{2k-1})=z\left(\dfrac{\alpha}{2}+\dfrac{1}{4}-g_{2k+1}\right)^2,
		\end{align*}
		and the result is immediate.
	\end{proof}	
	
	\begin{remark}
		Notice that for $\alpha= -1/2$ you get  (55) in \cite{DM22}. It is a  modified discrete Painlev\'e $ II $ equation (see \cite{Ramani,Magnus2}).
	\end{remark}
	
	\section{Holonomic Differential equation}\label{section5}
	\begin{teo}[Structure relation] The  polynomials $(P_n)_{n\geq 0}$ and $(Q_n)_{n\geq 0}$ orthogonal with respect to the linear functional $\elln$ and $x\elln$, respectively,   satisfy the differential-recurrence relation
		\begin{equation}\label{estructureeq}
			\begin{aligned}
				\phi(x)\partial_xP_{n+1}&=(n+1)P_{n+2}+\dfrac{1}{2}\left(a_{n+2}+a_{n+1}+b_{n+1}^2-(z+\alpha+2)b_{n+1}+z(1+\alpha)\right)P_{n+1}\\
				&+\left(b_{n+1}+b_n-(z+\alpha+2+n)\right)a_{n+1}P_{n}+a_{n+1}a_nP_{n-1}.
			\end{aligned}
		\end{equation}
		\begin{equation}\label{estructureeq1}
			\begin{aligned}
				\phi(x)\partial_xQ_{n+1}&=(n+1)Q_{n+2}+\dfrac{1}{2}\left(c_{n+2}+c_{n+1}+d_{n+1}^2-(z+\alpha+3)d_{n+1}+z(2+\alpha)\right)Q_{n+1} \\
				&+\left(d_{n+1}+d_n-(z+\alpha+3+n)\right)c_{n+1}Q_{n}+c_{n+1}c_nQ_{n-1},
			\end{aligned}
		\end{equation}
		where $\phi(x)=x^2-zx.$
	\end{teo}
	\begin{proof}
		Taking into account that $(P_k)_{k=0}^{n+2}$ is a basis for the linear sequence of polynomials of degree less than or equal to $n+2$, then
		$$\phi(x)\partial_xP_{n+1}=(n+1)P_{n+2}+\sum_{k=0}^{n+1}\alpha_{n+1,k}P_k.$$
		Using the orthogonality property and Pearson equation \eqref{pearson_eq}
		\begin{equation}\label{def.alpha}
			\begin{aligned}
				\alpha_{n+1,k}h_k&=\prodint{\elln,\phi(\partial_xP_{n+1})P_k}=\prodint{\elln,\phi\partial_x(P_{n+1}P_k)}-\prodint{\elln,\phi P_{n+1}(\partial_xP_k)}\\
				&=\prodint{\elln,\psi P_{n+1}P_k}-\prodint{\elln,\phi  P_{n+1}(\partial_xP_k)}.
			\end{aligned}
		\end{equation}
		This implies that $\alpha_{n+1,k}=0$ for all $0\leq k<n-1.$
		Now, from \eqref{eqq1}, \eqref{ppp} and \eqref{eqq2} we get
		\begin{equation*}
			\begin{aligned}
				\alpha_{n+1,n+1}h_{n+1}&=\prodint{\elln,\psi P_{n+1}P_{n+1}}-\prodint{\elln,\phi  P_{n+1}(\partial_xP_{n+1})}\\&=\dfrac{1}{2}\left(a_{n+2}+a_{n+1}+b_{n+1}^2-(z+\alpha+2)b_{n+1}+z(1+\alpha)\right)h_{n+1},\\
				\alpha_{n+1,n}h_n&=\prodint{\elln,\psi P_{n+1}P_n}-\prodint{\elln,\phi  P_{n+1}P^\prime_n}=\left(b_{n+1}+b_n-(z+\alpha+2+n)\right)a_{n+1}h_n.
			\end{aligned}
		\end{equation*}
		Finally,
		\begin{equation*}
			\begin{aligned}
				&\alpha_{n+1,n-1}h_{n-1}=\prodint{\elln,\psi P_{n+1}P_{n-1}}=\prodint{\elln,x^2 P_{n+1}P_{n-1}}=\prodint{\elln,xP_{n+1}P_{n}}=a_{n+1}h_n.
			\end{aligned}
		\end{equation*}
		Taking into account that $a_n=\dfrac{h_n}{h_{n-1}}$ we get \eqref{estructureeq}. \eqref{estructureeq1} follows in a similar way.
	\end{proof}
	\begin{remark}\label{remark1}
		From \eqref{11} and \eqref{cases} you get
		\begin{equation}\label{relambdartt}
			\begin{aligned}
				\lambda_{n,n-1}=-\sum_{k=0}^{n-1}b_k=
				-\dfrac{1}{2}\left(a_{n+1}+a_{n}+b_{n}^2-(z+\alpha+2n+2)b_{n}+z(\alpha+2n+1)\right),
			\end{aligned}
		\end{equation}
		Thus, we deduce
		\begin{equation*}
			\begin{aligned}
				\phi(x)\partial_xP_{n+1}&=(n+1)P_{n+2}+\left[(n+1)(b_{n+1}-z)+\sum_{k=0}^nb_k\right]P_{n+1} \\
				&+\left(b_{n+1}+b_n-(z+\alpha+2+n)\right)a_{n+1}P_{n}+a_{n+1}a_nP_{n-1}.
			\end{aligned}
		\end{equation*}
	\end{remark}
	\begin{coro}[Structure relation]\label{StrS}
Let $(S_n)_{n\geq 0}$ be the sequence of monic orthogonal polynomials with respect to the linear functional $\hn.$ Then the following structure relation holds
\begin{equation*}
\begin{aligned}
&\Phi(x;z)\partial_x(S_{n})=n\,S_{n+2}+\\
&\left[  \gamma_{n+2}\,\gamma_{n+1}+\gamma_{n}\,\gamma_{n-1}+(\gamma_{n+1}+\gamma_{n})^2-\left(z+\alpha+2\right)(\gamma_{n+1}+\gamma_{n})+z(1+\alpha)\right]S_{n}\\
&+\left(2\sum_{k=n-2}^{n+1}\gamma_k-(2z+2\alpha+n+2)\right)\gamma_{n}\,\gamma_{n-1}S_{n-2}+2(\gamma_{n}\,\gamma_{n-1}\,\gamma_{n-2}\,\gamma_{n-3})S_{n-4},
\end{aligned}
\end{equation*}
where $\Phi(x,z)=x(x^2-z).$
\end{coro}
\begin{proof}
Notice that from \eqref{estructureeq}
$$\phi(y)\partial_yP_{n+1}(y)=(n+1)P_{n+2}(y)+\sum_{k=n-1}^{n+1}\alpha_{n+1,k}P_k(y).$$
Suppose that $y=x^2$, then the above equation becomes
\begin{equation*}
\begin{aligned}
\dfrac{\phi(x^2)}{2x}\partial_x(P_{n+1}(x^2))&=(n+1)P_{n+2}(x^2)+\sum_{k=n-1}^{n+1}\alpha_{n+1,k}P_k(x^2)\\
\dfrac{\Phi(x)}{2}\partial_xS_{2n+2}(x)&=(n+1)S_{2n+4}(x)+\sum_{k=n-1}^{n+1}\alpha_{n+1,k}S_{2k}(x)
\end{aligned}
\end{equation*}
Now, taking into account the values of $\alpha_{n+1,k}$ in \eqref{def.alpha} and \eqref{relationsgamma} we get
\begin{equation*}
\begin{aligned}
&\Phi(x;z)\partial_x(S_{2n+2})=(2n+2)S_{2n+4}+\\
&\left[  \gamma_{2n+4}\,\gamma_{2n+3}+\gamma_{2n+2}\,\gamma_{2n+1}+(\gamma_{2n+3}+\gamma_{2n+2})^2-\left(z+\alpha+2\right)(\gamma_{2n+3}+\gamma_{2n+2})+z(1+\alpha)\right]S_{2n+2}\\
&+\left(2\sum_{k=2n}^{2n+3}\gamma_k-(2z+2\alpha+(2n+2)+2)\right)\gamma_{2n+2}\,\gamma_{2n+1}S_{2n}+2(\gamma_{2n+2}\,\gamma_{2n+1}\,\gamma_{2n}\,\gamma_{2n-1})S_{2n-2}.
\end{aligned}
\end{equation*}
Similarly, from \eqref{estructureeq1} we have
$$\phi(y)\partial_yQ_{n+1}(y)=(n+1)Q_{n+2}(y)+\sum_{k=n-1}^{n+1}\widetilde\alpha_{n+1,k}Q_k(y).$$
Suppose that $y=x^2$, then the above equation becomes
\begin{equation*}
\begin{aligned}
\dfrac{\phi(x^2)}{2x}\partial_x\left(xQ_{n+1}(x^2)\right)&=(n+1)xQ_{n+2}(x^2)+\sum_{k=n-1}^{n+1}\widetilde\alpha_{n+1,k}xQ_k(x^2)+\left(\dfrac{x^4-z}{2}\right)xQ_{n+1}(x^2),\\
\dfrac{\Phi(x)}{2}\partial_xS_{2n+3}(x)&=(n+1)S_{2n+5}(x)+\sum_{k=n-1}^{n+1}\widetilde\alpha_{n+1,k}S_{2k+1}(x)+\left(\dfrac{x^2-z}{2}\right)S_{2n+3}(x).
\end{aligned}
\end{equation*}
Since \eqref{x2S}, then taking into account the values of $\widetilde \alpha_{n+1,k}$ in \eqref{estructureeq1}, \eqref{relationsgamma} and shifting $n\to n-1$
\begin{equation*}
\begin{aligned}
&\Phi(x;z)\partial_x(S_{2n+1})=(2n+1)S_{2n+3}+\\
&\left[  \gamma_{2n+3}\,\gamma_{2n+2}+\gamma_{2n+1}\,\gamma_{2n}+(\gamma_{2n+2}+\gamma_{2n+1})^2-\left(z+\alpha+2\right)(\gamma_{2n+2}+\gamma_{2n+1})+z(1+\alpha)\right]S_{2n+1}\\
&+\left(2\sum_{k=2n-1}^{2n+2}\gamma_k-(2z+2\alpha+(2n+1)+2)\right)\gamma_{2n+1}\,\gamma_{2n}S_{2n-1}+2(\gamma_{2n+1}\,\gamma_{2n}\,\gamma_{2n-1}\,\gamma_{2n-2})S_{2n-3}
\end{aligned}
\end{equation*}
and the result follows.
\end{proof}
We will now obtain a lowering operator acting on the variable $x$ for $P_n(x,z).$
\begin{pro}\label{proL}
For $n\in\N_0$ let the operator $L_n$ be defined by
\begin{equation}\label{operatorL}
L_n=A_n(x;z)\partial_x-B_{n}(x;z),
\end{equation}
where
\begin{equation*}
\begin{aligned}
A_n(x;z)=\dfrac{\phi(x;z)}{C_n(x;z)},& &B_n(x;z)=\dfrac{\delta_n(x;z)}{C_{n}(x;z)}
\end{aligned}
\end{equation*}
where
$C_{n}(x;z)=a_{n+1}[x+b_{n+1}-(2n+\alpha+z+3)]$ and $\delta_n(x;z)=(n+1)(x-z)+\sum_{k=0}^nb_k-a_{n+1}.$  Then
$$L_nP_{n+1}=P_n,\qquad n\geq0.$$
\end{pro}
\begin{proof}
Using the structure relation \eqref{estructureeq} and \eqref{ttrr} we have
\begin{equation*}
\begin{aligned}
(\phi\partial_x-\alpha_{n+1,n+1})P_{n+1}&=(n+1)\left[(x-b_{n+1})P_{n+1}-a_{n+1}P_n\right]+\alpha_{n+1,n}P_n\\
&+\dfrac{\alpha_{n+1,n-1}}{a_n}\left[(x-b_n)P_n-P_{n+1}\right]\\
&=[(n+1)(x-b_{n+1})-a_{n+1}]P_{n+1}+[\alpha_{n+1,n}+(x-b_n-1)a_{n+1}]P_n.
\end{aligned}
\end{equation*}
So, remembering the values of $\alpha_{n+1,n+1}$ and $\alpha_{n+1,n}$ in \eqref{def.alpha} and taking into account \eqref{11} and Remark~\ref{remark1}, we get
\begin{equation*}
\begin{aligned}
&\left[\phi\partial_x-\dfrac{1}{2}\left[a_{n+2}+a_{n+1}+b_{n+1}^2-(z+\alpha+2n+2)b_{n+1}+z(1+\alpha)\right]-(n+1)x+a_{n+1}\right]P_{n+1}\\
&=a_{n+1}[x+b_{n+1}-(2n+\alpha+z+3)]P_n
\end{aligned}
\end{equation*}
and the result follows.
\end{proof}
\begin{pro}
For $n\in\N_0$, let $ D_{n+1}$ be the differential operator defined by
\begin{equation}\label{Dn}\begin{aligned}
D_{n+1}&=C_n\phi^2\partial_{xx}+\phi\left[(2x-z)C_n-a_{n+1}\phi-\delta_nC_n-C_n\left(\delta_{n-1}+\dfrac{C_{n-1}}{a_n}(x-b_n)\right)\right]\partial_{x}\\&+\left[\left(\delta_{n-1}+\dfrac{C_{n-1}}{a_n}(x-b_n)\right)\delta_nC_n+\dfrac{C_{n-1}}{a_n}C^2_{n}-(n+1)\phi C_{n}+a_{n+1}\phi\delta_n\right].
\end{aligned}\end{equation}
Then $D_{n+1}P_{n+1}=0.$
\end{pro}
\begin{proof}
From \eqref{operatorL} we get the following relation
\begin{equation}\label{edpL3}
\begin{aligned}
\phi  P^\prime_{n+1}&=\delta_nP_{n+1}+C_nP_{n}\\
&=\delta_nP_{n+1}+\dfrac{C_n}{a_{n+1}}\left[(x-b_{n+1})P_{n+1}-P_{n+2}\right]\\
&=\left(\delta_n+\dfrac{C_n}{a_{n+1}}(x-b_{n+1})\right)P_{n+1}-\dfrac{C_n}{a_{n+1}}P_{n+2}
\end{aligned}
\end{equation}
as well as
$$P_n=\dfrac{\phi }{C_n}P^\prime_{n+1}-\dfrac{\delta_n}{C_n}P_{n+1}.$$
Reversing \eqref{edpL3} and taking  the derivative
\begin{equation}\label{edpL}
\phi P^\prime_n=\phi\left[\left(\dfrac{\phi }{C_n}\right)^{\prime} P^{\prime}_{n+1}+\dfrac{\phi} {C_n}P^{\prime\prime}_{n+1}-\left(\dfrac{\delta_n}{C_n}\right)^{\prime}P_{n+1}-\dfrac{\delta_n}{C_n}P^{\prime}_{n+1}\right].
\end{equation}
On the other hand, making a shift $n+1\to n$ in \eqref{edpL} we get
\begin{equation}\label{edpL1}
\begin{aligned}
\phi  P^\prime_{n}&=\left(\delta_{n-1}+\dfrac{C_{n-1}}{a_{n}}(x-b_{n})\right)P_{n}-\dfrac{C_{n-1}}{a_{n}}P_{n+1}\\
&=\left(\delta_{n-1}+\dfrac{C_{n-1}}{a_{n}}(x-b_{n})\right)\left(\dfrac{\phi }{C_n}P^\prime_{n+1}-\dfrac{\delta_n}{C_n}P_{n+1}\right)-\dfrac{C_{n-1}}{a_{n}}P_{n+1}.
\end{aligned}
\end{equation}
Thus, if the equations \eqref{edpL} and \eqref{edpL1} are equated we get the result.
\end{proof}
\begin{remark}
Taking into account that the polynomials $(Q_{n})_{n\geq0}$ are orthogonal with respect to the linear functional $x\elln$, then the following fact are satisfied
\begin{enumerate}
\item 		If $\widehat{L}_n$ is the differential operator defined by
\begin{equation*}
\widehat{L}_n=\widehat{A}_n(x;z)\partial_x-\widehat{B}_{n}(x;z),
\end{equation*}
where
\begin{equation*}
\begin{aligned}
\widehat{A}_n(x;z)=\dfrac{\phi(x;z)}{\widehat{C}_n(x;z)},& 	&\widehat{B}_n(x;z)=\dfrac{\widehat{\delta}_n(x;z)}{\widehat{C}_{n}(x;z)}
\end{aligned}
\end{equation*}
with
$\widehat{C}_{n}(x;z)=c_{n+1}[x+d_{n+1}-(2n+\alpha+z+4)]$ and $\widehat{\delta}_n(x;z)=(n+1)(x-z)+\sum_{k=0}^nd_k-c_{n+1}.$  Then
$$\widehat{L}_nQ_{n+1}=Q_n,\qquad n\geq0.$$
\item If $ \widehat{D}_{n+1}$ is the differential operator defined by
\begin{equation*}
\begin{aligned}
\widehat{D}_{n+1}&=\widehat{C}_n\phi^2\partial^2_{x}+\phi\left[(2x-z)\widehat{C}_n-c_{n+1}\phi-\widehat{\delta}_n\widehat{C}_n-\widehat{C}_n\left(\widehat{\delta}_{n-1}+\dfrac{\widehat{C}_{n-1}}{c_n}(x-d_n)\right)\right]\partial_{x}\\
&+\left[\left(\widehat{\delta}_{n-1}+\dfrac{\widehat{C}_{n-1}}{c_n}(x-d_n)\right)\widehat{\delta}_n\widehat{C}_n+\dfrac{\widehat{C}_{n-1}}{c_n}\widehat{C}^2_{n}-(n+1)\phi \widehat{C}_{n}+c_{n+1}\phi\widehat{\delta}_n\right].
\end{aligned}\end{equation*}
Then $\widehat{D}_{n+1}Q_{n+1}=0.$
		\end{enumerate}
	\end{remark}
	\begin{pro}
		Let the operator $M_n$ be defined by
		\begin{equation*}
			M_n=W_n(x;z)\partial_x-V_{n}(x;z),\end{equation*}
		where $$W_n(x;z)=\dfrac{\Phi(x;z)}{Y_n(x;z)},\quad V_n(x;z)=\dfrac{\varphi_n(x;z)}{Y_n(x;z)}$$
		\begin{multline*}
			\varphi_n(x;z)=nx^2+\gamma_{n+2}\,\gamma_{n+1}-\gamma_{n}\,\gamma_{n-1}+(\gamma_{n+1}+\gamma_{n})^2\\ 			-\left(z+\alpha+2n+1\right)(\gamma_{n+1}+\gamma_{n})+z(1+\alpha)\end{multline*}
		and$$Y_n(x:z)=2\,\gamma_{n}\,\gamma_{n-1}(x^2+\gamma_n+\gamma_{n+1}-[z+\alpha+n+1]).$$
		Then
		$$M_nS_{n}=S_{n-2},\qquad n\geq 2.$$
	\end{pro}
	\begin{proof}
		The proof is a direct consequence of Corollary \ref{StrS} and \eqref{x2S}.
	\end{proof}
 	\section{The variable $z$}\label{section6}
	
	In this section we will study $\elln_{n},a_{n},b_{n}$ as functions of $z.$ From
	the definition of $\left\langle \elln  ,p\right\rangle ,$ we see that%
	\[
	\partial_{z}\left\langle \elln  ,p\right\rangle =\partial_{z}%
	{\displaystyle\int\limits_{0}^{z}}
	p\left(  x;z\right)  x^{\alpha}e^{-x}dx=%
	{\displaystyle\int\limits_{0}^{z}}
	\partial_{z}p\left(  x;z\right)  x^{\alpha}e^{-x}dx+p\left(  z;z\right)
	z^{\alpha}e^{-z},
	\]
	and hence%
	\begin{equation}
		\partial_{z}\left\langle \elln  ,p\right\rangle =\left\langle \elln  ,\partial
		_{z}p\right\rangle +p\left(  z;z\right)  z^{\alpha}e^{-z}. \label{Diego-Dz}%
	\end{equation}
	On the other hand, assuming that $\alpha>-1,$ we have%
	\begin{align*}
		\left\langle \elln  ,x\partial_{x}p\right\rangle  &  =%
		{\displaystyle\int\limits_{0}^{z}}
		x^{\alpha+1}e^{-x}\partial_{x}p\left(  x;z\right)  dx=p\left(  z;z\right)
		z^{\alpha+1}e^{-z}-%
		{\displaystyle\int\limits_{0}^{z}}
		p\left(  x;z\right)  \partial_{x}\left(  x^{\alpha+1}e^{-x}\right)  dx\\
		&  =p\left(  z;z\right)  z^{\alpha+1}e^{-z}-\left(  \alpha+1\right)
		{\displaystyle\int\limits_{0}^{z}}
		p\left(  x;z\right)  x^{\alpha}e^{-x}dx+%
		{\displaystyle\int\limits_{0}^{z}}
		p\left(  x;z\right)  x^{\alpha+1}e^{-x}dx,
	\end{align*}
	and therefore%
	\begin{equation}
		\left\langle \elln,x\partial_{x}p\right\rangle =p\left(  z;z\right)  z^{\alpha
			+1}e^{-z}-\left(  \alpha+1\right)  \left\langle \elln  ,p\right\rangle +\left\langle
		\elln,xp\right\rangle . \label{Diego-xDx}%
	\end{equation}
	Multiplying \eqref{Diego-Dz} by $z$ and comparing with \eqref{Diego-xDx}, we conclude that%
	\begin{equation}
		\vartheta\left\langle \elln  ,p\right\rangle =\left\langle \elln  ,\vartheta
		p\right\rangle +\left\langle \elln  ,x\partial_{x}p\right\rangle +\left(
		\alpha+1\right)  \left\langle \elln  ,p\right\rangle -\left\langle \elln  ,xp\right\rangle
		, \label{Diego-Dz Dx}%
	\end{equation}
	where $\vartheta$ denotes the differential operator%
	\begin{equation}
		\vartheta=z\partial_{z}. \label{Diego-theta}%
	\end{equation}

	\begin{pro}		The moments $\elln_{n}\left(  z\right) =\left\langle \elln,x^{n}\right\rangle $
		satisfy the differential-difference equation%
		\begin{equation}
			\vartheta\elln_{n}=\left(  n+\alpha+1\right) \elln_{n}-\elln_{n+1}, \label{Diego-l n+1}%
		\end{equation}
		and have the power series representation%
		\begin{equation}
			\elln_{n}\left(  z\right)  =z^{n+\alpha+1}%
			{\displaystyle\sum\limits_{k=0}^{\infty}}
			\dfrac{\left(  -1\right)  ^{k}}{n+\alpha+1+k}\dfrac{z^{k}}{k!}.
			\label{Diego-ln series}%
		\end{equation}
		
	\end{pro}
	
	\begin{proof}
		Using $p=x^{n}$ in \eqref{Diego-Dz Dx}, we obtain%
		\[
		\vartheta\elln_{n}=0+n \elln_{n}+\left(  \alpha+1\right) \elln_{n}-\elln_{n+1}.
		\]
On the other hand, if we replace
\[
\elln_{n}\left(  z\right)  =%
{\displaystyle\sum\limits_{k=0}^{\infty}}
c_{n,k}z^{n+\alpha+1+k}%
		\]
		in \eqref{Diego-l n+1}, we obtain%
		\[%
		{\displaystyle\sum\limits_{k=0}^{\infty}}
		c_{n+1,k}z^{n+\alpha+2+k}=\left(  n+\alpha+1\right)
		{\displaystyle\sum\limits_{k=0}^{\infty}}
		c_{n,k}z^{n+\alpha+1+k}-%
		{\displaystyle\sum\limits_{k=0}^{\infty}}
		\left(  n+\alpha+1+k\right)  c_{n,k}z^{n+\alpha+1+k},
		\]
		and therefore%
		\[
		c_{n+1,k-1}=\left(  n+\alpha+1\right)  c_{n,k}-\left(  n+\alpha+1+k\right)
		c_{n,k}=-kc_{n,k}.
		\]
		We conclude that%
		\[
		c_{n,k}=\dfrac{\left(  -1\right)  ^{k}}{\left(  n+\alpha+1+k\right)  k!}.
		\]
		
	\end{proof}
	
	If we define $\sigma_{n}\left(  z\right)  $ by
	\begin{equation}
		P_{n}\left(  x;z\right)  =x^{n}-\sigma_{n}\left(  z\right)  x^{n-1}+\mathcal{O}\left(
		x^{n-2}\right)  , \label{Diego-coeff}%
	\end{equation}
	then it follows from the three term recurrence relation \eqref{ttrrr} that%
	\begin{equation}
		b_{n}=\sigma_{n+1}-\sigma_{n}. \label{Diego-beta-sigma}%
	\end{equation}
	Note that in \eqref{rev} we wrote $\sigma_{n}=-\lambda_{n,n-1}$ and therefore  $nb_{n-1}+\lambda_{n,n-1}=\sigma_{n-1}-(n-1)b_{n-1}$.
	
	\begin{teo}
		Let $h_{n}$ be defined by $		h_{n}=\left\langle \elln  ,P_{n}^{2}\right\rangle , 
		$
		$a_{n},b_{n}$ be defined by \eqref{ttrrr}, and $\sigma_{n}$ be defined by
		\eqref{Diego-coeff}. We have%
		\begin{equation}
			\vartheta\sigma_{n}=\sigma_{n}-a_{n}, \label{Diego-sigma-gamma}%
		\end{equation}
		and%
		\begin{equation}
			\vartheta\ln\left(  h_{n}\right)  =2n+\alpha+1-b_{n}. \label{Diego-h-beta}%
		\end{equation}
		
	\end{teo}
	
	\begin{proof}
		In  \eqref{der} we shown that%
		\[
		\partial_{x}P_{n}=nP_{n-1}+\left[  \sigma_{n-1}-\left(  n-1\right)
		b_{n-1}\right]  P_{n-2}+\mathcal{O}\left(  x^{n-3}\right)  ,
		\]
		while \eqref{Diego-coeff} gives%
		\[
		\vartheta P_{n}=-\left(  \vartheta\sigma_{n}\right)  P_{n-1}++\mathcal{O}\left(
		x^{n-2}\right)  .
		\]
		Setting $p=P_{n-1}P_{n}$ in \eqref{Diego-Dz Dx}, we get%
		\[
		0=\vartheta\left\langle \elln  ,P_{n-1}P_{n}\right\rangle =\left\langle
		\elln,\vartheta\left(  P_{n-1}P_{n}\right)  \right\rangle +\left\langle
		\elln,x\partial_{x}\left(  P_{n-1}P_{n}\right)  \right\rangle -\left\langle
		\elln,xP_{n-1}P_{n}\right\rangle ,
		\]
		and since%
		\[
		\left\langle \elln  ,xP_{n-1}P_{n}\right\rangle =h_{n-1}a_{n},
		\]
		we obtain%
		\begin{align*}
			h_{n-1}a_{n}  &  =\left\langle \elln  ,P_{n}\vartheta P_{n-1}\right\rangle
			+\left\langle \elln  ,P_{n-1}\vartheta P_{n}\right\rangle +\left\langle
			\elln,xP_{n}\partial_{x}P_{n-1}\right\rangle +\left\langle \elln,xP_{n-1}\partial
			_{x}P_{n}\right\rangle \\
			&  =0-\left(  \vartheta\sigma_{n}\right)  h_{n-1}+0+b_{n-1}\left\langle
			\elln,P_{n-1}\partial_{x}P_{n}\right\rangle +a_{n-1}\left\langle \elln  ,P_{n-2}%
			\partial_{x}P_{n}\right\rangle \\
			&  =-\left(  \vartheta\sigma_{n}\right)  h_{n-1}+nb_{n-1}h_{n-1}%
			+a_{n-1}\left[  \sigma_{n-1}-\left(  n-1\right)  b_{n-1}\right]  h_{n-2}.
		\end{align*}
		Since
		\begin{equation}
			a_{n}=\dfrac{h_{n}}{h_{n-1}},\quad n\geq1, \label{Diego-gamma-h}%
		\end{equation}
		we conclude that%
		\[
		a_{n}=-\vartheta\sigma_{n}+nb_{n-1}+\sigma_{n-1}-\left(  n-1\right)  b_{n-1},
		\]
		or using \eqref{Diego-beta-sigma}%
		\[
		\vartheta\sigma_{n}=\sigma_{n}-a_{n}.
		\]

		Setting $p=P_{n}^{2}$ in \eqref{Diego-Dz Dx}, we get%
		\[
		\vartheta h_{n}=\vartheta\left\langle \elln,P_{n}^{2}\right\rangle =\left\langle
		\elln,\vartheta P_{n}^{2}\right\rangle +\left\langle \elln,x\partial_{x}P_{n}%
		^{2}\right\rangle +\left(  \alpha+1\right)  \left\langle \elln  ,P_{n}%
		^{2}\right\rangle -\left\langle \elln,xP_{n}^{2}\right\rangle ,
		\]
		and since%
		\[
		\left\langle \elln  ,xP_{n}^{2}\right\rangle =h_{n}b_{n},
		\]
		we obtain%
		\begin{align*}
			h_{n}b_{n}+\vartheta h_{n}  &  =2\left\langle \elln,P_{n}\vartheta P_{n}%
			\right\rangle +2\left\langle \elln,xP_{n}\partial_{x}P_{n}\right\rangle +\left(
			\alpha+1\right)  h_{n}\\
			&  =0+2na_{n}h_{n-1}+\left(  \alpha+1\right)  h_{n},
		\end{align*}
		and using \eqref{Diego-gamma-h}, we conclude that%
		\[
		\vartheta\ln\left(  h_{n}\right)  =2n+\alpha+1-b_{n}.
		\]
		
	\end{proof}
	
	As a consequence of the previous result, we see that $\left(  a_{n}%
	,b_{n}\right)  $ are solutions of the \emph{Toda lattice}.
	
	\begin{coro}
		The coefficients of the 3-term recurrence relation \eqref{ttrrr} are solutions
		of the differential-difference equations
		\begin{equation}
			\vartheta\ln\left(  a_{n}\right)  =2-\nabla b_{n},\quad\vartheta b_{n}%
			=b_{n}-\Delta a_{n}, \label{Diego-Toda}%
		\end{equation}
		with initial conditions
		\begin{equation}
			a_{0}=0,\quad b_{0}=\dfrac{\elln_{1}}{\elln_{0}}, \label{Diego-initial}%
		\end{equation}
		where the finite difference operators $\Delta,\nabla$ are defined by
		\begin{equation*}
			\Delta f\left(  n\right)  =f\left(  n+1\right)  -f\left(  n\right)
			,\quad\nabla f\left(  n\right)  =f\left(  n\right)  -f\left(  n-1\right)  .
		\end{equation*}
		
	\end{coro}
	
	\begin{proof}
		From \eqref{Diego-h-beta} and \eqref{Diego-gamma-h} we have%
		\begin{align*}
			\vartheta\ln\left(  a_{n}\right)   &  =\vartheta\ln\left(  h_{n}\right)
			-\vartheta\ln\left(  h_{n-1}\right)  =2n+\alpha+1-b_{n}-\left(  2n-1+\alpha
			-b_{n-1}\right) \\
			&  =2-\left(  b_{n}-b_{n-1}\right)  ,
		\end{align*}
		while \eqref{Diego-beta-sigma} and \eqref{Diego-sigma-gamma} give%
		\[
		\vartheta b_{n}=\vartheta\sigma_{n+1}-\vartheta\sigma_{n}=\sigma_{n+1}%
		-a_{n+1}-\left(  \sigma_{n}-a_{n}\right)  ,
		\]
		and using \eqref{Diego-beta-sigma} again we obtain%
		\[
		\vartheta b_{n}=b_{n}-\left(  a_{n+1}-a_{n}\right)  .
		\]
		
	\end{proof}
	
	\begin{teo}
		\label{Diego-Th1}Let the functions $a_{n}\left(  z\right)  ,b_{n}\left(  z\right)  $
		have the power series expansions
		\begin{equation}
			a_{n}\left(  z\right)  =%
			{\displaystyle\sum\limits_{k=0}^{\infty}}
			A_{n,k}z^{k},\quad b_{n}\left(  z\right)  =%
			{\displaystyle\sum\limits_{k=0}^{\infty}}
			B_{n,k}z^{k}. \label{Diego-a,b series}%
		\end{equation}
		Then, $A_{n,0}=0=B_{n,0},$%
		\begin{equation*}
			A_{n,1}=0,\quad B_{n,1}=\dfrac{2n\left(  n+\alpha+1\right)  +\alpha\left(
				\alpha+1\right)  }{\left(  2n+\alpha\right)  \left(  2n+\alpha+2\right)  },
		\end{equation*}%
		\begin{equation*}
			A_{n,2}=\dfrac{n^{2}(n+\alpha)^{2}}{(2n+\alpha-1)(2n+\alpha)^{2}(2n+\alpha
				+1)},\quad B_{n,2}=-\Delta A_{n,2}, 
		\end{equation*}
		and for all $k\geq3$%
		\begin{equation*}
			B_{n,k}=\dfrac{\Delta A_{n,k}}{1-k},\quad A_{n,k}=\dfrac{1}{2-k}%
			{\displaystyle\sum\limits_{j=1}^{k-2}}
			A_{n,k-j}\nabla B_{n,j}. 
		\end{equation*}
		
	\end{teo}
	
	\begin{proof}
		Using \eqref{Diego-ln series} in \eqref{Diego-initial}, we obtain $A_{0,k}=0$ for all
		$k\geq0$ and
		\begin{equation*}
			B_{0,0}=0,\quad B_{0,1}=\dfrac{\alpha+1}{\alpha+2},\quad B_{0,2}=-\dfrac
			{\alpha+1}{\left(  \alpha+2\right)  ^{2}\left(  \alpha+3\right)  }.
		\end{equation*}
		Using \eqref{Diego-a,b series} in \eqref{Diego-Toda}, we get%
		\[
		kB_{n,k}=B_{n,k}-\Delta A_{n,k},\quad kA_{n,k}=2A_{n,k}-%
		{\displaystyle\sum\limits_{j=0}^{k}}
		A_{n,k-j}\nabla B_{n,j},
		\]
		or%
		\begin{equation}
			\Delta A_{n,k}=\left(  1-k\right)  B_{n,k},\quad A_{n,k}=\dfrac{1}{2-k}%
			{\displaystyle\sum\limits_{j=0}^{k}}
			A_{n,k-j}\nabla B_{n,j}. \label{Diego-AB1}%
		\end{equation}
		Setting $k=0$ in \eqref{Diego-AB1}, we have%
		\[
		\Delta A_{n,0}=B_{n,0},\quad A_{n,0}=\dfrac{1}{2}A_{n,0}\nabla B_{n,0},
		\]
		and therefore%
		\[
		A_{n,0}=n\left(  n-1\right)  +c_{1}+c_{2}n,\quad B_{n,0}=2n+c_{2}.
		\]
		Since $A_{0,0}=0=B_{0,0},$ we conclude that%
		\begin{equation}
			A_{n,0}=0=B_{n,0},\quad n\geq0. \label{Diego-A0B0}%
		\end{equation}

		Using \eqref{Diego-A0B0} in \eqref{Diego-AB1}, we get%
		\begin{equation}
			\Delta A_{n,k}=\left(  1-k\right)  B_{n,k},\quad A_{n,k}=\dfrac{1}{2-k}%
			{\displaystyle\sum\limits_{j=1}^{k-1}}
			A_{n,k-j}\nabla B_{n,j}. \label{Diego-AB2}%
		\end{equation}
		Setting $k=1$ in \eqref{Diego-AB2} we have $\Delta A_{n,1}=0,$ and since
		$A_{0,1}=0,$ we see that
		\begin{equation*}
			A_{n,1}=0,\quad n\geq0. 
		\end{equation*}
		Thus, we can reduce \eqref{Diego-AB2} further and write%
		\begin{equation}
			\Delta A_{n,k}=\left(  1-k\right)  B_{n,k},\quad A_{n,k}=\dfrac{1}{2-k}%
			{\displaystyle\sum\limits_{j=1}^{k-2}}
			A_{n,k-j}\nabla B_{n,j}. \label{Diego-AB3}%
		\end{equation}
		Clearly, $k=2$ is a singular value for the system \eqref{Diego-AB3}, and therefore
		we will now use the Laguerre-Freud equations \eqref{111}--\eqref{112}.
		
		Setting
		\begin{equation}
			a_{n}\left(  z\right)  =%
			{\displaystyle\sum\limits_{k=2}^{\infty}}
			A_{n,k}z^{k},\quad b_{n}\left(  z\right)  =%
			{\displaystyle\sum\limits_{k=1}^{\infty}}
			B_{n,k}z^{k}. \label{Diego-ab series 1}%
		\end{equation}
		in \eqref{111}, we get%
		\[
		(2n+\alpha+2)B_{n,1}-(2n+\alpha-2)B_{n-1,1}=2,
		\]
		and since $B_{0,1}=\dfrac{\alpha+1}{\alpha+2},$ we conclude that%
		\begin{equation}
			B_{n,1}=\dfrac{2n\left(  n+\alpha+1\right)  +\alpha\left(  \alpha+1\right)
			}{\left(  2n+\alpha\right)  \left(  2n+\alpha+2\right)  }. \label{Diego-B1}%
		\end{equation}
		Using \eqref{Diego-ab series 1} in \eqref{112}, we get%
		\[
		(2n+\alpha+3)A_{n+1,2}-(2n+\alpha-1)A_{n,2}=(1-B_{n,1})B_{n,1},
		\]
		and since $A_{0,1}=0,$ we see that%
		\[
		A_{n,2}=\dfrac{1}{(2n+\alpha-1)(2n+\alpha+1)}%
		{\displaystyle\sum\limits_{k=0}^{n-1}}
		(2k+\alpha+1)(1-B_{k,1})B_{k,1},
		\]
		and using \eqref{Diego-B1} we obtain%
		\begin{equation}
			A_{n,2}=\dfrac{n^{2}(n+\alpha)^{2}}{(2n+\alpha-1)(2n+\alpha)^{2}(2n+\alpha+1)}.
			\label{Diego-A2}%
		\end{equation}

		Setting $k=2$ in \eqref{Diego-AB3} and using \eqref{Diego-A2}, we have
		\[
		B_{n,2}=-\Delta A_{n,2}=\dfrac{\widetilde{B}_{n,2}}{\left(  2n+\alpha-1\right)
			\left(  2n+\alpha\right)  ^{2}\left(  2n+\alpha+1\right)  \left(
			2n+\alpha+2\right)  ^{2}\left(  2n+\alpha+3\right)  },
		\]
		with%
		\begin{align*}
			\widetilde{B}_{n,2}  &  =-4\left(  2\alpha^{2}-1\right)  n^{4}-8\left(
			\alpha+1\right)  \left(  2\alpha^{2}-1\right)  n^{3}-2\left(  \allowbreak
			5\alpha^{4}+12\alpha^{3}+\alpha^{2}-6\alpha-2\right)  n^{2}\\
			&  -2\alpha\left(  \alpha+1\right)  \left(  \alpha^{3}+4\alpha^{2}%
			-\alpha-2\right)  n-\left(  \alpha-1\right)  \alpha^{2}\left(  \alpha
			+1\right)  ^{2}.
		\end{align*}
		Although proceeding this way one can find all coefficients $A_{n,k},B_{n,k},$
		the expressions become very big and difficult to write down.
	\end{proof}
	
	From \eqref{Diego-beta-sigma}, we have%
	\[
	\sigma_{n}\left(  z\right)  =%
	{\displaystyle\sum\limits_{k=0}^{n-1}}
	b_{k}\left(  z\right)  ,
	\]
	and using the results from the previous theorem we obtain%
	\begin{equation}
		\sigma_{n}\left(  z\right)  =\dfrac{n\left(  n+\alpha\right)  }{2n+\alpha
		}z-\dfrac{n^{2}(n+\alpha)^{2}}{(2n+\alpha-1)(2n+\alpha)^{2}(2n+\alpha+1)}%
		z^{2}+\mathcal{O}\left(  z^{3}\right)  \label{Diego-initial sigma}%
	\end{equation}
	as $z\rightarrow0.$ In the next theorem, we will find the higher order coefficients.
\begin{teo}
  If
		\begin{equation}
			\sigma_{n}\left(  z\right)  =%
			{\displaystyle\sum\limits_{k=0}^{\infty}}
			s_{n,k}z^{k}, \label{Diego-sigma series}%
		\end{equation}
		then%
		\begin{equation}
			s_{n,0}=0,\quad s_{n,1}=\dfrac{n\left(  n+\alpha\right)  }{2n+\alpha},\quad
			s_{n,2}=-\dfrac{n^{2}(n+\alpha)^{2}}{(2n+\alpha-1)(2n+\alpha)^{2}(2n+\alpha
				+1)}, \label{Diego-sn first}%
		\end{equation}
		and for $k\geq3$
		\begin{equation}
			s_{n,k}=\dfrac{-1}{\left(  k-1\right)  \left(  k-2\right)  }%
			{\displaystyle\sum\limits_{j=1}^{k-2}}
			\left(  k-1-j\right)  s_{n,k-j}\Delta\nabla s_{n,j}. \label{Diego-snk req}%
		\end{equation}
		
	\end{teo}
	
	\begin{proof}
		From \eqref{Diego-beta-sigma}, \eqref{Diego-sigma-gamma}, and \eqref{Diego-Toda}, we have%
		\[
		\vartheta\ln\left(  \sigma_{n}-\vartheta\sigma_{n}\right)  =\vartheta
		\ln\left(  a_{n}\right)  =2-\nabla b_{n}=2-\Delta\nabla\sigma_{n},
		\]
		and using \eqref{Diego-theta}, we get
		\begin{equation}
			z^{2}\sigma_{n}^{\prime\prime}=\left(  z\sigma_{n}^{\prime}-\sigma_{n}\right)
			\left(  2-\Delta\nabla\sigma_{n}\right)  . \label{Diego-sigmaeq}%
		\end{equation}
		Using \eqref{Diego-sigma series} in \eqref{Diego-sigmaeq}, we obtain%
		\[
		k\left(  k-1\right)  s_{n,k}=2\left(  k-1\right)  s_{n,k}-%
		{\displaystyle\sum\limits_{j=0}^{k}}
		\left(  k-j-1\right)  s_{n,k-j}\Delta\nabla s_{n,j},
		\]
		or%
		\begin{equation}
			\left(  k-1\right)  \left(  k-2\right)  s_{n,k}=-%
			{\displaystyle\sum\limits_{j=0}^{k}}
			\left(  k-j-1\right)  s_{n,k-j}\Delta\nabla s_{n,j}. \label{Diego-snj1}%
		\end{equation}

		From \eqref{Diego-initial sigma} we get \eqref{Diego-sn first}, and since $s_{n,0}=0$ we
		can reduce \eqref{Diego-snj1} to
		\[
		\left(  k-1\right)  \left(  k-2\right)  s_{n,k}=-%
		{\displaystyle\sum\limits_{j=1}^{k-2}}
		\left(  k-j-1\right)  s_{n,k-j}\Delta\nabla s_{n,j}%
		\]
		and \eqref{Diego-snk req} follows.
	\end{proof}
	
	It is well known (\cite{Chihara}) that if we define the Hankel determinants
	$\mathcal{H}_{n}\left(  z\right)  $ by $\mathcal{H}_{0}=1$ and
	\begin{equation*}
		\mathcal{H}_{n}\left(  z\right)  =\underset{0\leq i,j\leq n-1}{\det}\left(
		\mu_{i+j}\right)  ,\quad n\geq1, 
	\end{equation*}
	then,%
	\begin{equation*}
		\mathcal{H}_{n}=%
		{\displaystyle\prod\limits_{k=0}^{n-1}}
		h_{k},\quad n\geq1. 
	\end{equation*}
	From \eqref{Diego-beta-sigma} and \eqref{Diego-h-beta}, we get%
	\[
	\sigma_{n+1}-\sigma_{n}=b_{n}=2n+\alpha+1-\vartheta\ln\left(  h_{n}\right)
	=2n+\alpha+1-\vartheta\ln\left(  \mathcal{H}_{n+1}\right)  +\vartheta
	\ln\left(  \mathcal{H}_{n}\right)  ,
	\]
	and since $\sigma_{0}=0=\vartheta\ln\left(  \mathcal{H}_{0}\right)  ,$ we
	conclude that%
	\begin{equation}
		\sigma_{n}=n\left(  n+\alpha\right)  -\vartheta\ln\left(  \mathcal{H}%
		_{n}\right)  =-\vartheta\ln\left(  z^{-n\left(  n+\alpha\right)  }%
		\mathcal{H}_{n}\right)  . \label{Diego-sigma Hn}%
	\end{equation}

	In \cite{MR4093808}, the author studied power series expansions of Hankel
	determinants of the form \eqref{Diego-Hn} with $\mu_{n}\in\mathbb{F}.$ Using the
	results obtained in that paper, one can show that
	\begin{equation}
		\mathcal{H}_{n}\left(  z\right)  =g_{n}\left(  z\right)  \left(
		n+\alpha\right)  ^{n}z^{n\left(  n+\alpha\right)  }%
		{\displaystyle\prod\limits_{j=0}^{n-1}}
		\dfrac{\left(  j!\right)  ^{2}}{\left(  j+\alpha+1\right)  ^{j+1}\left(
			n+j+\alpha\right)  ^{n-j}},\quad n\geq1, \label{Diego-Hn}%
	\end{equation}
	with%
	\[
	g_{1}\left(  z\right)  =\left(  \alpha+1\right)  z^{-\left(  \alpha+1\right)
	}\elln_{0}\left(  z\right)  =\left(  \alpha+1\right)
	{\displaystyle\sum\limits_{k=0}^{\infty}}
	\dfrac{\left(  -1\right)  ^{k}}{\alpha+1+k}\dfrac{z^{k}}{k!},
	\]
	and in general
	\begin{equation}
		g_{n}\left(  z\right)  =1-\dfrac{n\left(  n+\alpha\right)  }{2n+\alpha}%
		z+\dfrac{n^{2}\left(  n+\alpha\right)  ^{2}}{\left(  2n+\alpha-1\right)
			\left(  2n+\alpha+1\right)  }\dfrac{z^{2}}{2}+\mathcal{O}\left(  z^{3}\right)  ,
		\label{Diego-gn series}%
	\end{equation}
	as $z\rightarrow0.$ If we use \eqref{Diego-gn series} in \eqref{Diego-sigma Hn}, we
	recover \eqref{Diego-sn first}.
	
	\section{Asymptotic analysis}\label{section7}
	
	Using \eqref{Diego-sn first} and \eqref{Diego-snk req} (we omit writing the HUGE explicit
	expressions for the higher order coefficients), we obtain%
	\begin{equation}%
		\begin{tabular}[c]{l}%
			$s_{n,1}=\dfrac{n}{2}+\dfrac{\alpha}{4}-\dfrac{\alpha^{2}}{8n}+\dfrac{\alpha^{3}%
			}{16n^{2}}-\dfrac{\alpha^{4}}{32n^{3}}+\mathcal{O}\left(  n^{-4}\right)  ,$\\[10 pt]
			$s_{n,2}=-\dfrac{1}{16}+\dfrac{2\alpha^{2}-1}{64n^{2}}+\dfrac{\alpha\left(
				1-2\alpha^{2}\right)  }{64n^{3}}+\mathcal{O}\left(  n^{-4}\right)  ,$\\[10 pt]
			$s_{n,3}=-\dfrac{\alpha^{2}}{128n^{3}}+\dfrac{3\alpha^{3}}{256n^{4}}%
			-\dfrac{\alpha^{2}\left(  4\alpha^{2}+5\right)  }{512n^{5}}+\mathcal{O}\left(n^{-6}\right)  ,$\\[10 pt]
			$s_{n,4}=\dfrac{2\alpha^{2}-1}{1024n^{4}}+\dfrac{\alpha\left(  1-2\alpha Z	^{2}\right)  }{512n^{5}}+\mathcal{O}\left(  n^{-6}\right)  ,$%
		\end{tabular}
		\label{Diego-snk first large}%
	\end{equation}
	as $n\rightarrow\infty.$ As it was observed in \cite{MR4136730}, the
	coefficients $s_{n,k}$ (after some $k)$ form an \emph{asymptotic sequence} as
	$n\rightarrow\infty.$
	
	\begin{pro}
		Let the coefficients $s_{n,k}$ be defined by \eqref{Diego-sn first} and
		\eqref{Diego-snk req}. Then, for all $k\geq3$%
		\[
		s_{n,k}\sim\dfrac{\left(  -1\right)  ^{k}\left(  4\alpha^{2}-1\right)
			-1}{2^{2k+3}}n^{-k},\quad n\rightarrow\infty.
		\]
		
	\end{pro}
	
	\begin{proof}
		Suppose that for all $k\geq3$
		\begin{equation}
			s_{n,k}=c_{k}n^{\theta_{k}}+\mathcal{O}\left(  n^{\theta_{k}-1}\right)  ,\quad
			n\rightarrow\infty, \label{Diego-snk large}%
		\end{equation}
		with $\theta_{k}>\theta_{k+1}.$ Using the binomial theorem, we have%
		\begin{equation}
			\Delta\nabla s_{n,k}=\theta_{k}\left(  \theta_{k}-1\right)  c_{k}n^{\theta
				_{k}-2}+\mathcal{O}\left(  n^{\theta_{k}-3}\right)  ,\quad n\rightarrow\infty.
			\label{Diego-diff snk}%
		\end{equation}
		If $s_{n,k}$ is an asymptotic sequence (in $k),$ we expect that the lower
		terms will have the higher order (in $n).$ Therefore, the dominant terms in
		\eqref{Diego-snk req} should be the boundary terms%
		\begin{equation}
			-\left(  k-1\right)  \left(  k-2\right)  s_{n,k}\sim\left(  k-2\right)
			s_{n,k-1}\Delta\nabla s_{n,1}+s_{n,2}\Delta\nabla s_{n,k-2},\quad
			n\rightarrow\infty. \label{Diego-req  snk asy}%
		\end{equation}

		Using \eqref{Diego-snk large} and \eqref{Diego-diff snk} in \eqref{Diego-req snk asy}, we get
		(to leading order)
		\begin{equation}
			-\left(  k-1\right)  \left(  k-2\right)  c_{k}n^{\theta_{k}}\sim-\dfrac
			{\alpha^{2}}{4}\left(  k-2\right)  c_{k-1}n^{\theta_{k-1}-3}-\dfrac{1}%
			{16}\theta_{k-2}\left(  \theta_{k-2}-1\right)  c_{k-2}n^{\theta_{k-2}-2}
			\label{Diego-req asy 1}%
		\end{equation}
		since from \eqref{Diego-snk first large} we know that%
		\[
		\Delta\nabla s_{n,1}\sim-\dfrac{\alpha^{2}}{4n^{3}},\quad s_{n,2}\sim-\dfrac
		{1}{16}+\dfrac{2\alpha^{2}-1}{64n^{2}},\quad n\rightarrow\infty.
		\]
		From \eqref{Diego-snk first large} we have
		\[
		\theta_{k-1}-3=-7<-5=\theta_{k-2}-2,\quad k=5,
		\]
		so we assume that
		\begin{equation}
			\theta_{k-1}-3<\theta_{k-2}-2\text{ \ for all \ }k\geq3, \label{Diego-cond theta}%
		\end{equation}
		and the leading terms in \eqref{Diego-req asy 1} give
		\begin{equation}
			\theta_{k}=\theta_{k-2}-2, \label{Diego-theta req}%
		\end{equation}
		and%
		\begin{equation}
			-\left(  k-1\right)  \left(  k-2\right)  c_{k}=-\dfrac{1}{16}\theta
			_{k-2}\left(  \theta_{k-2}-1\right)  c_{k-2}. \label{Diego-ck req}%
		\end{equation}

		Solving \eqref{Diego-theta req} with initial conditions $\theta_{3}=-3,$ $\theta
		_{4}=-4,$ we obtain
		\begin{equation}
			\theta_{k}=-k,\quad k\geq3, \label{Diego-thetak}%
		\end{equation}
		which is consistent with \eqref{Diego-cond theta}. Using \eqref{Diego-thetak} in
		\eqref{Diego-ck req}, we get
		\[
		c_{k}=\dfrac{1}{16}c_{k-2},\quad k\geq5,
		\]
		and since from \eqref{Diego-snk first large} we have the initial conditions
		\[
		c_{3}=-\dfrac{\alpha^{2}}{128},\quad c_{4}=\dfrac{2\alpha^{2}-1}{1024},
		\]
		we conclude that%
		\[
		c_{k}=\dfrac{\left(  -1\right)  ^{k}\left(  4\alpha^{2}-1\right)  -1}{2^{2k+3}%
		},\quad k\geq3.
		\]
		Note that%
		\[
		c_{2k-1}=-\dfrac{\alpha^{2}}{2^{4k-1}},\quad c_{2k}=\dfrac{2\alpha^{2}%
			-1}{2^{4k+2}},\quad k\geq2.
		\]
		
	\end{proof}
	
	As a consequence of the previous result, we can use \eqref{Diego-snk first large}
	and write the first terms in the asymptotic expansion of $\sigma_{n}\left(
	z\right)  $ as $n\rightarrow\infty$%
	
	\begin{equation*}
		\sigma_{n}\left(  z\right)  =\dfrac{z}{2}n+\dfrac{\left(  4\alpha-z\right)
			z}{16}-\dfrac{\alpha^{2}z}{8}n^{-1}+z\dfrac{2\alpha^{2}(z+2\alpha)-z}{64}%
		n^{-2}+\mathcal{O}\left(  n^{-3}\right)  . 
	\end{equation*}
	From \eqref{Diego-beta-sigma}, \eqref{Diego-sigma-gamma}, and \eqref{Diego-a,b series}, we have%
	\[
	A_{n,k}=\left(  1-k\right)  s_{n,k},\quad B_{n,k}=\Delta s_{n,k},
	\]
	and therefore%
	\begin{equation}
		a_{n}\left(  z\right)  =\dfrac{z^{2}}{16}+\dfrac{\left(  1-2\alpha^{2}\right)
			z^{2}}{64}n^{-2}+\alpha z^{2}\dfrac{\alpha(z+2\alpha)-1}{64}n^{-3}+\mathcal{O}\left(
		n^{-4}\right)  , \label{Diego-a asy}%
	\end{equation}%
	\begin{equation}
		b_{n}\left(  z\right)  =\dfrac{z}{2}-\dfrac{\alpha^{2}z}{8}n^{-2}+z\dfrac
		{z-2\alpha^{2}(z+2\alpha+2)}{32}n^{-3}+\mathcal{O}\left(  n^{-4}\right)  , \label{Diego-b asy}%
	\end{equation}
	as $n\rightarrow\infty$, in agreement with Theorem \ref{Diego-Th1}.
	
	Using \eqref{Diego-a asy} and \eqref{Diego-b asy} in the 3-term recurrence \eqref{ttrrr},
	we see that (to leading order)%
	\[
	xP_{n}\sim P_{n+1}+\dfrac{z}{2}P_{n}+\dfrac{z^{2}}{16}P_{n-1},\quad
	n\rightarrow\infty
	\]
	and therefore%
	\[
	P_{n}\left(  x;z\right)  \sim\left(  \dfrac{z}{4}\right)  ^{n}U_{n}\left(
	\dfrac{2}{z}x-1\right)  ,\quad n\rightarrow\infty
	\]
	where $U_{n}\left(  x\right)  $ denotes the Chebyshev polynomial of the second
	kind, defined by $U_{0}=1,$ $U_{1}=2x$ and
	\[
	U_{n+1}=2xU_{n}-U_{n-1},\quad n\geq1.
	\]
	In particular, it follows that if we denote by $x_{n,k}$ the zeros of
	$P_{n}\left(  x;z\right)  ,$ then
	\[
	x_{n,k}\sim\dfrac{z}{2}\left[  \cos\left(  \dfrac{k\pi}{n+1}\right)  +1\right]
	,\quad1\leq k\leq n,\quad n\rightarrow\infty.
	\]
 \section{Painlev\'{e} equations}\label{section8}
	
	In this section we will derive a nonlinear differential equation satisfied by
	the function $\sigma_{n}\left(  z\right)$ defined in \eqref{Diego-coeff} To begin with, we will rewrite
	the Laguerre-Freud equations \eqref{111}--\eqref{112} in a factorized form.
	
	\begin{pro}
		\label{Diego-New LF}Let $\mathfrak{L}_{0},\mathfrak{L}_{1},\mathfrak{L}_{2},\mathfrak{L}_{3}:\mathbb{N}_{0}\rightarrow\mathbb{R}\left[  \left[
		z\right]  \right]  $ be defined by%
		\begin{align}
		\mathfrak{L}_{0}\left(  n\right)&=\left(  b_{n}-z\right)  b_{n}+\left(
		2n+3+\alpha+z-b_{n}-b_{n+1}\right)  a_{n+1}-\left(  2n-1+\alpha+z-b_{n-1}%
		-b_{n}\right)  a_{n}, \notag\\
			\mathfrak{L}_{1}\left(  n\right)&=a_{n}\left(  \alpha+2n-1-b_{n-1}\right)
			\left(  \alpha+1+2n-b_{n}\right)  -\left(  a_{n}-\sigma_{n}\right)
			^{2},\label{Diego-LF1}\\
			\mathfrak{L}_{2}\left(  n\right)  &=\left(  \alpha+1+2n\right)  z+a_{n}%
			+a_{n+1}-\left(  \alpha+2+2n+z-b_{n}\right)  b_{n}-2\sigma_{n},\label{Diego-LF2}\\
			\mathfrak{L}_{3}\left(  n\right)  &=n\left(  \alpha+n\right)  z-a_{n}\left(\alpha+2n+z-b_{n-1}-b_{n}\right)  -\left(  \alpha+2n\right)  \sigma
			_{n}.\label{Diego-LF3}%
		\end{align}
		Then,%
		\[
		\mathfrak{L}_{0}\left(  n\right)  =\mathfrak{L}_{1}\left(  n\right)
		=\mathfrak{L}_{2}\left(  n\right)  =\mathfrak{L}_{3}\left(  n\right)  =0,\quad
		n\geq0.
		\]
		
	\end{pro}
	
	\begin{proof}
		\begin{enumerate}
			\item The identity $\mathfrak{L}_{0}\left(  n\right)  =0$ for all $n\geq0$ is just the Laguerre-Freud equation \eqref{112}. The proof
			that $\mathfrak{L}_{1}\left(  n\right)  =0$ for all $n\geq0$
			follows immediately from Theorem~\ref{Teo4.4}.
			\item  We have
			\begin{align*}
				\mathfrak{L}_{2}\left(  n\right)  -\mathfrak{L}_{2}\left(  n-1\right)   &
				=\allowbreak2z+a_{n+1}-a_{n-1}-b_{n}\left(  2n+2+\alpha+z-b_{n}\right) \\
				&  +b_{n-1}\left(  2n+\alpha+z-b_{n-1}\right)  -2\left(  \sigma_{n}-\sigma_{n-1}\right)  ,
			\end{align*}
			and using \eqref{Diego-beta-sigma} we obtain the Laguerre-Freud equation \eqref{111}%
			\[
			\allowbreak2z+a_{n+1}-a_{n-1}-b_{n}\left(  2n+2+\alpha+z-b_{n}\right)+b_{n-1}\left(  2n+\alpha+z-b_{n-1}-2\right)  =0.\]
			Therefore, $\mathfrak{L}_{2}\left(  n\right)  $ must be constant and setting
			$n=0,$ we get%
			\[\mathfrak{L}_{1}\left(  0\right)  =\left(  \alpha+1\right)  z+a_{0}%
			+a_{1}-\left(  \alpha+2+z\right)  b_{0}+b_{0}^{2}-2\sigma_{0},\]
			or, since $a_{0}=0=\sigma_{0},$%
			\[
			\mathfrak{L}_{2}\left(  0\right)  =\left(  \alpha+1\right)  z+a_{1}-\left(
			\alpha+2+z\right)  b_{0}+b_{0}^{2}.
			\]
			From the three-term recurrence relation \eqref{ttrrr} we see that $P_{1}\left(
			x\right)  =x-b_{0}$, and using \eqref{recmoments}, \eqref{Diego-gamma-h}, and
			\eqref{Diego-initial}, we obtain%
			\begin{align*}
				a_{1}  &  =\dfrac{h_{1}}{h_{0}}=\dfrac{1}{\elln_{0}}\prodint{\elln,\left(  x-b_{0}\right)^{2}}  =\dfrac{\elln_{2}-2b_{0}\elln_{1}+b_{0}^{2}\elln_{0}}{\elln_{0}}=\dfrac{\elln_{2}%
				}{\elln_{0}}-b_{0}^{2}\\
				&  =\dfrac{\left(  z+2+\alpha\right) \elln_{1}-z\left(  \alpha+1\right) \elln_{0}}{\elln_{0}}-b_{0}^{2}=\left(  z+2+\alpha\right)  b_{0}-z\left(  \alpha+1\right)-b_{0}^{2}.
			\end{align*}
			Thus, $\mathfrak{L}_{2}\left(  0\right)  =0$ and we conclude that
			$\mathfrak{L}_{2}\left(  n\right)  =0$ for all $n\geq0.$
			\item We have%
			\begin{gather*}
				\mathfrak{L}_{3}\left(  n+1\right)  -\mathfrak{L}_{3}\left(  n\right)
				=z\left(  2n+\alpha+1\right)  -a_{n+1}\left(  \alpha+2n+2+z-b_{n}%
				-b_{n+1}\right)  \\
				+a_{n}\left(\alpha+2n+z-b_{n-1}-b_{n}\right)  -\left(  \alpha+2n+2\right)
				\sigma_{n+1}+\left(  \alpha+2n\right)  \sigma_{n}.
			\end{gather*}
			Using \eqref{Diego-beta-sigma} we obtain%
			\begin{gather*}
				\mathfrak{L}_{3}\left(  n+1\right)  -\mathfrak{L}_{3}\left(  n\right)
				=z\left(  2n+\alpha+1\right)  -a_{n+1}\left(  \alpha+2n+2+z-b_{n}%
				-b_{n+1}\right)  \\
				+a_{n}\left(  \alpha+2n+z-b_{n-1}-b_{n}\right)  -\left(  2n+\alpha+2\right)
				b_{n}-2\sigma_{n}.
			\end{gather*}
			On the other hand,%
			\begin{align*}
				0  & =\mathfrak{L}_{2}\left(  n\right)  -\mathfrak{L}_{0}\left(  n\right)
				=\left(  \alpha+1+2n\right)  z-a_{n+1}\left(  2n+2+z+\alpha-b_{n}%
				-b_{n+1}\right)  \\
				& +a_{n}\left(2n+z+\alpha-b_{n}-b_{n-1}\right)  -\left(  2n+\alpha+2\right)
				b_{n}-2\sigma_{n},
			\end{align*}
			and therefore $\mathfrak{L}_{3}\left(  n\right)  $ is constant. Since
			$\mathfrak{L}_{3}\left(  0\right)  =0,$ we conclude that $\mathfrak{L}_{3}\left(  n\right)=0$ for all $n\geq0.$
		\end{enumerate}
	\end{proof}
	Motivated by the work of Mourad E. H.  Ismail and Yang Chen on  \emph{ladder operators}
	\cite{MR1616931}, we define the functions (see also \cite{MR2525311, MR2570057, MR2736877, MR2265830, MR2581382,MR2676782})
	\begin{equation}
		R_{n}\left(  z\right)  =\dfrac{b_{n}\left(  z\right)  -\left(  2n+\alpha
			+1\right)  }{z},\quad r_{n}\left(  z\right)  =\dfrac{a_{n}\left(  z\right)
			-\sigma_{n}\left(  z\right)  }{z}.\label{Diego-Rr}%
	\end{equation}
Note that in terms of $R_{n}\left(  z\right)  ,r_{n}\left(  z\right)  $ we can write Proposition \ref{proL}
	\[\partial_{x}P_{n}\left(  x;z\right)  +\widetilde{B}_{n}\left(  x;z\right)
	P_{n}\left(  x;z\right)  =a_{n}\left(  z\right)  \widetilde{A}_{n}\left(
	x;z\right)  P_{n-1}\left(  x;z\right),
	\]
	where%
	\[
	\widetilde{A}_{n}=\dfrac{R_{n}}{x-z}+\dfrac{1-R_{n}}{x},\quad\widetilde{B}%
	_{n}=\dfrac{r_{n}}{x-z}-\dfrac{r_{n}+n}{x}.
	\]
	Therefore,  $R_{n}\left(  z\right)  ,r_{n}\left(  z\right)  $ are residues
	of the lowering operator $a_{n}\left(  z\right)  \widetilde{A}_{n}\left(
	x;z\right)  L_{n-1}$ at the pole $x=z.$ Using $R_{n}\left(  z\right)
	,r_{n}\left(  z\right)  $ we can rewrite the results of Proposition
	\ref{Diego-New LF}.
	\begin{pro}
		For all $n\geq0$ we have
		\begin{equation}
			r_{n}^{2}=a_{n}R_{n}R_{n-1},\label{Diego-Rr1}%
		\end{equation}%
		\begin{equation}
			r_{n+1}+r_{n}=-\left(2n+1+\alpha-z+zR_{n}\right)  R_{n},\label{Diego-Rr2}%
		\end{equation}
		and
		\begin{equation}
			\left(  r_{n}+n\right)  \left(  r_{n}+n+\alpha\right)  =a_{n}\left(
			R_{n}-1\right)  \left(  R_{n-1}-1\right)  .\label{Diego-Rr3}%
		\end{equation}
	\end{pro}
	
	\begin{proof}
		Using \eqref{Diego-Rr} in \eqref{Diego-LF1} and \eqref{Diego-LF2}, we get%
		\[\mathfrak{L}_{1}\left(  n\right)  =a_{n}R_{n}R_{n-1}-r_{n}^{2},\]
		and%
		\[\mathfrak{L}_{2}\left(  n\right)  =\left(  \alpha+1+2n\right)  z+z\left(
		r_{n}+r_{n+1}\right)  -\left(  z-zR_{n}+1\right)  b_{n}+\sigma_{n+1}%
		+\sigma_{n}-2\sigma_{n},\]
		or using \eqref{Diego-beta-sigma}%
		\[	\mathfrak{L}_{2}\left(  n\right)  =r_{n+1}+r_{n}+\left(  2n+1+\alpha
		-z+zR_{n}\right)  R_{n}.\]
		Finally, using \eqref{Diego-Rr} in \eqref{Diego-LF3}, we have%
		\[\mathfrak{L}_{3}\left(  n\right)  =n\left(  \alpha+n\right)  z-a_{n}\left(
		z-2n-\alpha-zR_{n}-zR_{n-1}\right)  -\left(  \alpha+2n\right)  \left(a_{n}-zr_{n}\right),
		\]
		$\allowbreak$and therefore%
		\begin{equation}
			\alpha\left(  n+r_{n}\right)  +n^{2}+2nr_{n}=a_{n}\left(  1-R_{n-1}-R_{n}\allowbreak\right)  .\label{Diego-Rr4}%
		\end{equation}
		Adding \eqref{Diego-Rr1} and \eqref{Diego-Rr4}, we obtain%
		\[\alpha\left(  n+r_{n}\right)  +n^{2}+2nr_{n}+r_{n}^{2}=a_{n}\left(
		1-R_{n-1}-R_{n}\allowbreak\right)  +a_{n}R_{n}R_{n-1},\]
		and \eqref{Diego-Rr3} follows.
	\end{proof}
	As a corollary, we obtain an expression for $a_{n}$ in terms of $R_{n}$ and
	$r_{n}.$
	\begin{coro}
		We have%
		\begin{equation}
			a_{n}=\dfrac{r_{n}^{2}}{R_{n}}-\dfrac{\left(  r_{n}+n\right)\left(r_{n}+n+\alpha\right)  }{R_{n}-1},\quad n\geq0.\label{Diego-arR}%
		\end{equation}
		
	\end{coro}
	
	\begin{proof}
		From \eqref{Diego-Rr1} and \eqref{Diego-Rr3} we get%
		\[
		a_{n}=\dfrac{\left(  r_{n}+n\right)  \left(  r_{n}+n+\alpha\right)  }{\left(
			R_{n}-1\right)  \left(  R_{n-1}-1\right)  }=\dfrac{\left(  r_{n}+n\right)
			\left(  r_{n}+n+\alpha\right)  }{\left(  R_{n}-1\right)  \left(  \dfrac
			{r_{n}^{2}}{a_{n}R_{n}}-1\right)  },
		\]
		and, therefore, %
		\[
		a_{n}=\dfrac{r_{n}^{2}+R_{n}\left(  n\alpha+2nr_{n}+\alpha r_{n}+n^{2}\right)
		}{R_{n}\left(  1-R_{n}\right)  },
		\]
		from which \eqref{Diego-arR} follows.
	\end{proof}
	
	In terms of $R_{n}$ and $r_{n},$ the Toda equations \eqref{Diego-Toda} take the very
	simple form%
	\begin{equation}
		\partial_{z}\ln\left(  a_{n}\right)  =-\nabla R_{n},\quad\partial_{z}%
		b_{n}=-\Delta r_{n},\label{Diego-Rr Toda}%
	\end{equation}
	since%
	\[
	2-\nabla b_{n}=2-\left(  zR_{n}+2n+\alpha+1\right)  +zR_{n-1}+2n+\alpha
	-1=z\left(  R_{n-1}-R_{n}\right)  ,
	\]
	and%
	\[
	b_{n}-\Delta a_{n}=b_{n}-\left(  zr_{n+1}+\sigma_{n+1}\right)  +zr_{n}%
	+\sigma_{n}=z\left(  r_{n}-r_{n+1}\right)  .
	\]

	Using our previous results, we can write a system of first order ODEs for
	$R_{n}$ and $r_{n}.$
	
	\begin{teo}
		We have%
		\begin{equation}
			\vartheta R_{n}=2r_{n}+\left(  2n+\alpha-z+zR_{n}\right)  R_{n},\label{Diego-DR}%
		\end{equation}
		and%
		\begin{equation}
			\vartheta r_{n}=n\left(  n+\alpha\right)  +\left(  2n+\alpha\right)
			r_{n}+\dfrac{r_{n}^{2}}{R_{n}}+\allowbreak\dfrac{\left(  n+r_{n}\right)  \left(
				n+\alpha+r_{n}\right)  }{R_{n}-1}.\label{Diego-Dr}%
		\end{equation}
		
	\end{teo}
	
	\begin{proof}
		According to \eqref{Diego-Rr} and \eqref{Diego-Rr Toda}, we get%
		\[
		r_{n}-r_{n+1}=\partial_{z}\left(  zR_{n}+2n+\alpha+1\right)  =R_{n}%
		+zR_{n}^{\prime},
		\]
		and from \eqref{Diego-Rr2} we see that%
		\[
		-r_{n+1}=r_{n}+\left(  2n+1+\alpha-z+zR_{n}\right)  R_{n},
		\]
		and \eqref{Diego-DR} follows.
		
		The first equation in \eqref{Diego-Rr Toda} and \eqref{Diego-arR} gives%
		\begin{gather*}
			R_{n-1}-R_{n}=\partial_{z}\ln\left[  \dfrac{r_{n}^{2}+\left(  n\alpha
				+2nr_{n}+\alpha r_{n}+n^{2}\right)  R_{n}}{R_{n}\left(  1-R_{n}\right)
			}\right]  \\
			=\dfrac{\partial_{z}\left[  r_{n}^{2}+\left(  n\alpha+2nr_{n}+\alpha
				r_{n}+n^{2}\right)  R_{n}\right]  }{r_{n}^{2}+\left(  n\alpha+2nr_{n}+\alpha
				r_{n}+n^{2}\right)  R_{n}}-\dfrac{\partial_{z}\left[  R_{n}\left(
				1-R_{n}\right)  \right]  }{R_{n}\left(  1-R_{n}\right)  }\\
			=\dfrac{2r_{n}r_{n}^{\prime}+\left(  n\alpha+2nr_{n}+\alpha r_{n}+n^{2}\right)
				R_{n}^{\prime}+R_{n}\left(  2n+\alpha\right)  r_{n}^{\prime}}{r_{n}%
				^{2}+\left(  n\alpha+2nr_{n}+\alpha r_{n}+n^{2}\right)  R_{n}}-\dfrac{\left(
				1-2R_{n}\right)  R_{n}^{\prime}}{R_{n}\left(  1-R_{n}\right)  }.
		\end{gather*}
		Multiplying by $z$ and using \eqref{Diego-DR}, we obtain%
		\begin{gather*}
			z\left(  R_{n-1}-R_{n}\right)  =\dfrac{\left[  2r_{n}+\left(  2n+\alpha\right)
				R_{n}\right]  \vartheta r_{n}}{r_{n}^{2}+\left(  n\alpha+2nr_{n}+\alpha
				r_{n}+n^{2}\right)  R_{n}}\\
			+\left[  \dfrac{n\alpha+2nr_{n}+\alpha r_{n}+n^{2}}{r_{n}^{2}+\left(
				n\alpha+2nr_{n}+\alpha r_{n}+n^{2}\right)  R_{n}}+\dfrac{1-2R_{n}}{R_{n}\left(
				R_{n}-1\right)  }\right]  \left[  2r_{n}+\left(  2n+\alpha-z+zR_{n}\right)
			R_{n}\right]  .
		\end{gather*}
		$\allowbreak$Using \eqref{Diego-Rr1}, we have
		\begin{gather*}
			z\left(  \dfrac{r_{n}^{2}}{a_{n}R_{n}}-R_{n}\right)  =\dfrac{\left[
				2r_{n}+\left(  2n+\alpha\right)  R_{n}\right]  \vartheta r_{n}}{r_{n}%
				^{2}+\left(  n\alpha+2nr_{n}+\alpha r_{n}+n^{2}\right)  R_{n}}\\
			+\left[  \dfrac{n\alpha+2nr_{n}+\alpha r_{n}+n^{2}}{r_{n}^{2}+\left(
				n\alpha+2nr_{n}+\alpha r_{n}+n^{2}\right)  R_{n}}+\dfrac{1-2R_{n}}{R_{n}\left(
				R_{n}-1\right)  }\right]  \left[  2r_{n}+\left(  2n+\alpha-z+zR_{n}\right)
			R_{n}\right].
		\end{gather*}
		From \eqref{Diego-arR}, we get%
		\begin{gather*}
			z\left(  \dfrac{r_{n}^{2}\left(  1-R_{n}\right)  }{r_{n}^{2}+R_{n}\left(
				n\alpha+2nr_{n}+\alpha r_{n}+n^{2}\right)  }-R_{n}\right)  =\dfrac{\left[
				2r_{n}+\left(  2n+\alpha\right)  R_{n}\right]  \vartheta r_{n}}{r_{n}%
				^{2}+\left(  n\alpha+2nr_{n}+\alpha r_{n}+n^{2}\right)  R_{n}}\\
			+\left[  \dfrac{n\alpha+2nr_{n}+\alpha r_{n}+n^{2}}{r_{n}^{2}+\left(
				n\alpha+2nr_{n}+\alpha r_{n}+n^{2}\right)  R_{n}}+\dfrac{1-2R_{n}}{R_{n}\left(
				R_{n}-1\right)  }\right]  \left[  2r_{n}+\left(  2n+\alpha-z+zR_{n}\right)
			R_{n}\right]  .
		\end{gather*}
		Solving for $\vartheta r_{n},$ we obtain%
		\[
		\vartheta r_{n}=\dfrac{2R_{n}r_{n}^{2}+n^{2}R_{n}^{2}-r_{n}^{2}+2nR_{n}%
			^{2}r_{n}+\alpha R_{n}^{2}r_{n}+n\alpha R_{n}^{2}}{R_{n}\left(  R_{n}%
			-1\right)  },
		\]
		and the result follows.
	\end{proof}
	
	We have now all the necessary ingredients to find a nonlinear ODE satisfied by
	$\sigma_{n}\left(  z\right)  .$
	
	\begin{teo}
		The function $\sigma_{n}\left(  z\right)  $ is a solution of the nonlinear ODE%
		\begin{equation}
			\left(  zy^{\prime\prime}\right)  ^{2}+4\left(  y^{\prime}\right)
			^{2}(n-y^{\prime})(n+\alpha-y^{\prime})-\left[  n\left(  n+\alpha\right)
			-y+\left(  2y^{\prime}+z-2n-\alpha\right)  y^{\prime}\right]  ^{2}%
			=0,\label{Diego-ODE}%
		\end{equation}
		where $y^{\prime}=\partial_{z}y.$
	\end{teo}
	
	\begin{proof}
		Using \eqref{Diego-sigma-gamma} and \eqref{Diego-Rr}, we have%
		\[
		z\sigma_{n}^{\prime}=\sigma_{n}-a_{n}=-zr_{n},
		\]
		and, therefore,  $\sigma_{n}^{\prime}=-r_{n}.$ Thus, we can rewrite \eqref{Diego-ODE}
		in terms of $r_{n}\left(  z\right) .$  In such a sense we need to show that%
		\begin{equation}
			\left(  \vartheta r_{n}\right)  ^{2}+4r_{n}^{2}(n+r_{n})(n+\alpha
			+r_{n})-\left[  n\left(  n+\alpha\right)  -a_{n}+\left(  2n+\alpha
			+2r_{n}\right)  r_{n}\right]  ^{2}=0.\label{Diego-ODE1}%
		\end{equation}
		From \eqref{Diego-arR} and \eqref{Diego-Dr}, we get%
		\begin{align*}
	& \vartheta r_{n}+n\left(  n+\alpha\right)  -a_{n}+\left(  2n+\alpha+2r_{n}\right) r_{n}\\
& =n\left(  n+\alpha\right)  +\left(  2n+\alpha\right)  r_{n}+\dfrac{r_{n}^{2}}{R_{n}}+\allowbreak\dfrac{\left(  n+r_{n}\right)\left(n+\alpha+r_{n}\right)  }{R_{n}-1}\\
& +n\left(  n+\alpha\right)  -\left[  \dfrac{r_{n}^{2}}{R_{n}}-\dfrac{\left(
r_{n}+n\right)  \left(  r_{n}+n+\alpha\right)  }{R_{n}-1}\right]  +\left(2n+\alpha+2r_{n}\right)  r_{n}\\
& =\dfrac{2\left(  n+r_{n}\right)  \left(  n+\alpha+r_{n}\right)  R_{n}}{R_{n}-1},
\end{align*}
and%
\begin{align*}
& \vartheta r_{n}-n\left(  n+\alpha\right)  +a_{n}-\left(  2n+\alpha
+2r_{n}\right)  r_{n}\\
& =n\left(  n+\alpha\right)  +\left(  2n+\alpha\right)  r_{n}+\dfrac{r_{n}^{2}}{R_{n}}+\allowbreak\dfrac{\left(  n+r_{n}\right)  \left(  n+\alpha
+r_{n}\right)  }{R_{n}-1}\\
& -n\left(  n+\alpha\right)  +\dfrac{r_{n}^{2}}{R_{n}}-\dfrac{\left(
r_{n}+n\right)  \left(  r_{n}+n+\alpha\right)}{R_{n}-1}-\left(
2n+\alpha+2r_{n}\right)  r_{n}\\
& =-\dfrac{2r_{n}^{2}\left(  R_{n}-1\right)  }{R_{n}}.
\end{align*}
We conclude that%
\begin{align*}
\left(  \vartheta r_{n}\right)  ^{2}-\left[  n\left(  n+\alpha\right)
-a_{n}+\left(  2n+\alpha+2r_{n}\right)  r_{n}\right]  ^{2}
& =-\dfrac{2r_{n}^{2}\left(  R_{n}-1\right)  }{R_{n}}\dfrac{2\left(
n+r_{n}\right)  \left(  n+\alpha+r_{n}\right)  R_{n}}{R_{n}-1}\\
& =-4r_{n}^{2}\left(  n+r_{n}\right)  \left(n+\alpha+r_{n}\right),
\end{align*}
and \eqref{Diego-ODE1} follows.
\end{proof}
We can rewrite \eqref{Diego-ODE} as%
\begin{equation}
\left(  zY^{\prime\prime}\right)  ^{2}-\left(  Y+Y^{\prime}\left(  2Y^{\prime}-z\right)\right)^{2}+4%
{\displaystyle\prod\limits_{k=1}^{4}}\left(  Y^{\prime}+v_{k}\right)  =0,\label{Diego-ODE2}%
\end{equation}
$\allowbreak$with%
\[	v_{1}=-\dfrac{1}{4}\left(  2n+\alpha\right)  =v_{2},\quad v_{3}=\dfrac{1}%
{4}\left(  2n-\alpha\right) ,\quad v_{4}=\dfrac{1}{4}\left(  2n+3\alpha\right),\]
and%
\begin{align*}
Y  & =-\sigma_{n}+\left(  v_{3}-v_{1}\right)  \left(  v_{4}-v_{1}\right)
-v_{1}z-2v_{1}^{2}\\
& =-\sigma_{n}+\dfrac{1}{2}n\left(  n+\alpha\right)  -\dfrac{\alpha^{2}}%
{8}+\dfrac{1}{4}\left(  2n+\alpha\right)  z.
\end{align*}
The nonlinear ODE \eqref{Diego-ODE2} is the Jimbo-Miwa-Okamoto form of Painlev\'{e}
	V \cite{MR1885665,MR914314}.
\section{Zeros}\label{section9}
It is well known that the zeros of orthogonal polynomials with respect to a positive definite linear functional are real, simple, and located in the interior of the convex hull of the support \cite{Chi,GMM21}.
With this in mind, let $\{x_{n,k}(z)\}_{k=1}^{n}$   be the 	the zeros of $P_{n}(x)$  in an increasing order, i. e.,
\begin{equation}\label{zerosP_n}
P_{n}(x_{n,k}(z),z)=0
\end{equation}
with
$$x_{n,1}(z)<x_{n,2}(z)<\cdots<x_{n,n}(z).$$

Next we will focus our attention on an electrostatic interpretation of the zeros of the polynomial $P_{n}(x ;z)$ in terms of the energy associated with a logarithmic potential and next we will study the dynamics of them in terms of the parameter $z.$
\subsection{Electrostatic interpretation}

Evaluating the operator $D_{n+1}$ defined in \eqref{Dn} at $x = x_{n+1,k},$ we get
\begin{multline*}	\left.\dfrac{\partial_x^2P_{n+1}}{\partial_xP_{n+1}}\right|_{x=x_{n+1,k}}=-\dfrac{(2x_{n+1,k}-z)}{\phi(x_{n+1,k})}+\dfrac{a_{n+1}}{C_n(x_{n+1,k})}\\+\dfrac{(x_{n,k}-b_n)C_{n-1}(x_{n,k};z)}{a_{n}\phi(x_{n,k};z)} +\dfrac{\delta_n(x_{n,k};z)+\delta_{n-1}(x_{n,k};z)}{\phi(x_{n,k};z)},
\end{multline*}
where $C_{n}(x;z)$ and $\delta_n(x;z)$ were defined in Proposition \ref{proL}. Taking into account \eqref{relambdartt}  the above expression reads
\begin{equation*}
\begin{aligned}
\left.\dfrac{\partial_x^2P_{n+1}}{\partial_xP_{n+1}}\right|_{x=x_{n+1,k}}
&=-\dfrac{1}{x_{n+1,k}-z}-\dfrac{1}{x_{n+1,k}}+\dfrac{1}{x_{n+1,k}-\beta_n}+1-\dfrac{\alpha}{x_{n+1,k}},
\end{aligned}
\end{equation*}
where $\beta_n=-b_{n+1}+(2n+\alpha+z+3)$.
Observe that the weight function associated with $\elln$ is $w(x)=x^{\alpha}e^{-x}=e^{-x+\alpha\ln x}.$ If we define $v(x)=x-\alpha\ln x,$ (the external potential \cite[Section~3.5]{Mourad}), then
\begin{equation}\label{equielect1}
\begin{aligned}
\left.\dfrac{\partial_x^2P_{n+1}}{\partial_xP_{n+1}}\right|_{x=x_{n+1,k}}=-\dfrac{1}{x_{n+1,k}-z}-\dfrac{1}{x_{n+1,k}}+\dfrac{1}{x_{n+1,k}-\beta_n}+v^\prime (x_{n+1,k}).
\end{aligned}
\end{equation}
\begin{remark}
If $\{\widetilde{x}_{n+1,k}(z)\}_{k=1}^{n+1}$   are  the zeros of $Q_{n+1}(x)$  in an increasing order, i.e. $\widetilde{x}_{n+1,1}<\widetilde{x}_{n+1,2}<\cdots<\widetilde{x}_{n+1,n+1}$, then we have
\begin{equation}\label{equielectQ1}
\begin{aligned}
\left.\dfrac{\partial_x^2Q_{n+1}}{\partial_x Q_{n+1}}\right|_{x=\widetilde{x}_{n+1,k}}=-\dfrac{1}{\widetilde{x}_{n+1,k}-z}-\dfrac{1}{\widetilde{x}_{n+1,k}}+\dfrac{1}{\widetilde{x}_{n+1,k}-\widetilde{\beta}_n}+1-\dfrac{\alpha+1}{\widetilde{x}_{n+1,k}},
\end{aligned}
\end{equation}
where $\widetilde{\beta}_n=-d_{n+1}+(2n+\alpha+z+4)$.
\end{remark}
\begin{teo}
The zeros of $P_{n+1}(x; z)$ are the equilibrium points of $n+1$ unit charged particles located in the interval $(0, z)$ under the influence of the  external potential
$$	V_P(x)=\ln\dfrac{1}{|x-z|}+\ln\dfrac{1}{|x|}-\ln\dfrac{1}{|x-\beta_n|}+v(x), \quad 0<x<z.$$
\end{teo}
\begin{proof} If $P_{n+1}(x)=\prod_{k=1}^{n+1}(x-x_{n+1,k})$
then \cite[Chapter 10]{GMM21} we have
\begin{equation*}
\sum_{\substack{j=1 \\ j\ne k} }^{n+1}\dfrac{2}{x_{n+1,j}-x_{n+1,k}}=-\dfrac{\partial^2_xP_{n+1}(x_{n+1,k})}{\partial_xP_{n+1} (x_{n+1,k})}.
\end{equation*}
Taking into account  \eqref{equielect1}
\begin{equation*}
\sum_{\substack{j=1 \\ j\ne k} }^{n+1}\dfrac{2}{x_{n+1,j}-x_{n+1,k}}-\dfrac{1}{x_{n+1,k}-z}-\dfrac{1}{x_{n+1,k}}+\dfrac{1}{x_{n+1,k}-\beta_n}+v^\prime (x_{n+1,k})=0.
\end{equation*}
It is not difficult to see that the above system of equations can be written as
$$\dfrac{\partial E}{\partial \,x_{k}}=0, \quad 1\le k \le n+1,$$
where
\begin{align*}
E(x_{1},\ldots, &x_{n+1})=  2\sum_{1\le j < k \le n+1}\ln\dfrac{1}{|x_{j}-x_{k}|}+\\
&\sum_{k=1}^{n+1}\left[\ln\dfrac{1}{|x_{k}-z|}+\ln\dfrac{1}{|x_{k}|}-\ln\dfrac{1}{|x_{k}-\beta_n|}+v(x_{k})\right].
\end{align*}
\end{proof}
Notice that it involves the  ``natural" potential $v$ associated with the weight functions, positive unit charges located at the ends of the support of the weight and a charge at the extra point $\beta_{n}$ which is greater than $z$ for all $n\geq0$ (see Remark \ref{reamrkbeta}).\\
Let $\{y_{n+1,k}(z)\}_{k=1}^{n+1}$  be the
the zeros of $S_{n+1}(x)$ in an increasing order, i.e. $y_{n+1,1}<y_{n+1,2}<\cdots<y_{n+1,n+1}.$ Notice that if $x_{n+1,k}$ is a zero of $P_{n+1}(x)$, then $\sqrt{x_{n+1,k}}$ and $-\sqrt{x_{n+1,k}}$ are zeros of $S_{2n+2}(x)$, and if $\widetilde{x}_{n+1,k}$ is a zero of $Q_{n+1}(x)$, then $\sqrt{\widetilde{x}_{n+1,k}}$ and $-\sqrt{\widetilde{x}_{n+1,k}}$ are zeros of $S_{2n+1}(x)$, Note also that $0$ is  a zero of $S_{2n+1}(x)$.  Moreover if we take $y=x^2$, then
$$\begin{aligned}
\partial_{y}P_{n+1}(y)=\dfrac{1}{2x}\partial_x\left(P_{n+1}(x^2)\right)\quad &\text{ and }\quad  \partial^2_{y}P_{n+1}(y)=\dfrac{1}{4x^2}\partial^2_x\left(P_{n+1}(x^2)\right)-\dfrac{1}{4x^3}\partial_x\left(P_{n+1}(x^2)\right)\\
\partial_{y}P_{n+1}(y)=\dfrac{1}{2x}\partial_x\left(S_{2n+2}(x)\right)\quad &\text{ and }\quad  \partial^2_{y}P_{n+1}(y)=\dfrac{1}{4x^2}\partial^2_x\left(S_{2n+2}(x)\right)-\dfrac{1}{4x^3}\partial_x\left(S_{2n+2}(x)\right)
\end{aligned},$$
Thus, if $y_{2n+2}$ is a zero of $S_{2n+2}$, then taking into account the above we get
\begin{multline*}
\dfrac{1}{2y_{2n+2,k}}\left.\dfrac{\partial_x^2S_{2n+2}}{\partial_xS_{2n+2}}\right|_{x=y_{2n+2,k}}-\dfrac{1}{2y_{{2n+2,k}}^2}=-\dfrac{1}{y^2_{2n+2,k}-z}-\dfrac{1}{y^2_{2n+2,k}}+\dfrac{1}{y^2_{2n+2,k}-\beta_n}+1-\dfrac{\alpha}{y^2_{2n+2,k}},
\end{multline*}
or, equivalently
\begin{multline}\label{electrorhopar}
\left.\dfrac{\partial_x^2S_{2n+2}}{\partial_xS_{2n+2}}\right|_{x=y_{2n+2,k}}=\dfrac{1}{y_{2n+2,k}-\sqrt{\beta_n}}+\dfrac{1}{y_{2n+2,k}+\sqrt{\beta_n}}\\-\dfrac{1}{(y_{2n+2,k}-\sqrt{z})}-\dfrac{1}{(y_{2n+2,k}+\sqrt{z})}+2y_{2n+2,k}-\dfrac{2\alpha+1}{y_{2n+2,k}}.
\end{multline}
On the other hand, if  we take $y=x^2,$ then from  \eqref{equielectQ1}
$$\partial_{y}Q_{n+1}(y)=\dfrac{1}{2x}\partial_x\left(S_{2n+3}(x)\right)\quad \text{ and }\quad  \partial^2_{y}Q_{n+1}(y)=\dfrac{1}{4x^2}\partial^2_x\left(S_{2n+3}(x)\right)-\dfrac{1}{4x^3}\partial_x\left(S_{2n+3}(x)\right)$$
With this in mind, we have that
\begin{equation*}
\begin{aligned}
\partial_x S_{2n+3}(x)&=Q_{n+1}(y)+2y\,\partial_y Q_{n+1}(y)\\
x\partial^2_x S_{2n+3}(x)&=3\,\partial_yQ_{n+1}(y)+2y\,\partial^2_y Q_{n+1}(y).
\end{aligned}
\end{equation*}
Thus,
\begin{multline*}
\dfrac{1}{2x}\left.\dfrac{\partial_x^2S_{2n+3}(x)}{\partial_xS_{2n+3}(x)}\right|_{x=y_{2n+3,k}}=\dfrac{3}{2y}+\left.\dfrac{\partial_y^2Q_{n+1}(y)}{\partial_yQ_{n+1}(y)}\right|_{y=y^2_{2n+3,k}}\\
\\=\dfrac{3}{2y^2_{2n+3,k}}-\dfrac{1}{y^2_{2n+3,k}-z}-\dfrac{1}{y^2_{2n+3,k}}+\dfrac{1}{y^2_{2n+3,k}-\widetilde{\beta}_n}+1-\dfrac{\alpha+1}{y^2_{2n+3,k}},
\end{multline*}
or, equivalently
\begin{multline}\label{electrorhopar1}
\left.\dfrac{\partial_x^2S_{2n+3}(x)}{\partial_xS_{2n+3}(x)}\right|_{x=y_{2n+3,k}}=\dfrac{1}{y_{2n+3,k}-\sqrt{\widetilde{\beta}_n}}+\dfrac{1}{y_{2n+3,k}+\sqrt{\widetilde{\beta}_n}}\\-\dfrac{1}{(y_{2n+3,k}-\sqrt{z})}-\dfrac{1}{(y_{2n+3,k}+\sqrt{z})}+2y_{2n+3,k}-\dfrac{2\alpha+1}{y_{2n+3,k}}.
\end{multline}
Therefore we get the following result
\begin{teo}
The zeros of $S_{n+1}(x; z)$ are located at the equilibrium points of $n+1$ unit charged particles located in the interval $(-\sqrt{z}, \sqrt{z})$ under the influence of the potential
\begin{equation*}
V_S(x)=\ln\dfrac{1}{|x-\sqrt{z}|}+\ln\dfrac{1}{|x+\sqrt{z}|}-\ln\dfrac{1}{|x-\sqrt{\rho_n}|}-\ln\dfrac{1}{|x+\sqrt{\rho_n}|}+\widehat{v}(x),\quad -\sqrt{z}<x<\sqrt{z},\end{equation*}
where $\widehat{v}(x)=x^2-(2\alpha+1)\ln|x|$ is the  ``natural" potential associated with the weight function, two positive unit charges are  located at the ends of the support of the weight and two charges at the extra points $\pm \sqrt{\rho_{n}}.$
\end{teo}
\begin{remark}
Notice that we recover Theorem~21 in \cite{DM22}.
\end{remark}
Taking into account that  $\rho_n=\sqrt{(n+z+\alpha+1)-\gamma_{n}-\gamma_{n+1}}$ and the notation in Corollary \ref{coro_simetricos_painleve}  we get $\rho^2-z=g_n+g_{n+1}.$ Replacing  the above in \eqref{painlevesimetricos} we get
\begin{equation}\label{inequation}
\gamma_n(z-\rho_{n}^2)(z-\rho_{n-1}^2)>0,\qquad n\geq1.
\end{equation}
Since by definition $\gamma_0=0,$ then $z-\rho_0^2=\gamma_1-(\alpha+1)$. On the other hand, from definition
$$\gamma_1=\dfrac{1}{\mathfrak{h}_0}\int_{-\sqrt{z}}^{\sqrt{z}}|x|^{2n+3}e^{-x^2}dx=-\dfrac{z^{\alpha+1}e^{-z}}{\mathfrak{h}_0}+\alpha+1.$$ Therefore,
$$z-\rho_0^2=-\dfrac{z^{\alpha+1}e^{-z}}{\mathfrak{h}_0}<0.$$
Finally, from \eqref{inequation} we conclude that 
$\rho_n$ is outside the interval $(-\sqrt{z},\sqrt{z})$ for all $n\geq0.$
\begin{remark}\label{reamrkbeta}
Notice that from \eqref{electrorhopar} and \eqref{electrorhopar1} we get
$$\rho^2_{2n}=\beta_n \quad \text{and}\quad \rho^2_{2n+1}=\widetilde{\beta}_n, \quad n\geq0.$$  From here
$$z-\beta_n=z-\rho^2_{2n}<0 \quad \text{and}\quad z-\widetilde{\beta}_n=z-\rho^2_{2n+1}<0.$$
Therefore, $\beta_n>z$ and $\widetilde{\beta}_n>z$.
\end{remark}
\medskip
\subsection{Dynamical behavior of zeros}

Next,  we are interested to study the  motion of zeros of the time dependent   polynomials $(P_{n}(x,z))_{n\geq 0}$, which are orthogonal with respect to the linear functional $\elln$ defined in \eqref{laguerrefunc}. We follow the approach in \cite{Ismail-Wen11}.

 \begin{pro}[\cite{Ismail-Wen11}] \label{Prodiff Q and P}
Let  $ (P_n)_{n\geq 0}$ be the sequence of orthogonal  polynomials with respect to $\elln$. Then, for every $n\geq 1$,
\begin{equation*}
z\dot{P}_n(x,z)=a_n(t)P_{n-1}(x,z)-x\,\partial_xP_{n}(x,z)+nP_n(x,z).
\end{equation*}
where  the dot `\.{}' means derivative with respect to  the  variable $z$. 	
\end{pro}
\begin{proof}
If we define the linear functional  $\un$  by
	\begin{equation*}
		\langle\un,p(x)\rangle=\int_0^{1} p(x)\,x^{\alpha} e^{-zx}dx, \quad p(x)\in\mathbb{P},  \quad \alpha>-1,\quad z>0.
	\end{equation*}
	It is not difficult to check that if $(Q_{n}(x,z))_{n\geq 0}$ is the sequence of monic orthogonal polynomials with respect $\un$, where \begin{equation}\label{rel. QP}
	Q_{n}(x,z)=\dfrac{1}{z^n}P_{n}(zx,z),\quad n\geq 0.
	\end{equation}
	Moreover
$$zQ_{n}(x,z)=Q_{n+1}(x,z)+\dfrac{b_n(z)}{z}Q_{n}(x,z)+\dfrac{a_n(z)}{z^2}Q_{n-1}(x,z).$$
Taking into account that the measure associated with $\un$ is supported on $[0,1]$ and can be written as $d\mu=e^{-zx}d\nu(x)$ with $d\nu(x)=x^\alpha dx$, then    $(Q_n(x,z))_{n\geq 0}$ satisfies the following differential property (see \cite{Ismail-Wen11})
\begin{equation*}
\dot{Q}_n(x,z)=\dfrac{a_n(z)}{z^2}Q_{n-1}(x,z). 
		\end{equation*}
  Thus, taking into account \eqref{rel. QP}
  we get
$$x \partial_y\left. P_n(y,z)\right|_{y=zx}+\dot{P}_n(zx,z)-\dfrac{n}{z}P_{n}(zx,x)=\dfrac{a_n(z)}{z}P_{n-1}(zx,z).$$
Finally, multiplying both sides of the equation by $t$ and changing $zx\to x$ we get the result.
  \end{proof}	
\begin{coro} From Proposition \ref{Prodiff Q and P}, we get that  the differential operator
\begin{equation}\label{opediff2}
\begin{aligned}	
\mathcal{L}_2=\dfrac{1}{a_n(z)}\left(z\partial_z+x\,\partial_x-n\right)
\end{aligned}
\end{equation}
is an annihilation operator for the sequence  $(P_n)_{n\geq 0}$.
\end{coro}
Now, we are going to provide a differential equation which is satisfied by the zeros of $P_{n}(z,x)$.
Differentiating \eqref{zerosP_n}  with respect to $z$
$$\partial_xP_n(x,z)\left|_{x=x_{n,k}(z)}\right.\partial_z x_{n,k}(z)\,+\dot {P}_n(x,z)\left|_{x=x_{n,k}(z)}\right.=0.$$
The above implies that (see \cite{Ismail-Wen11} for a more general case)
\begin{equation}\label{diffy}
\dot{x}_{n,k}(z)=-\dfrac{\dot {P}_n(x,z)\left|_{x=x_{n,k}(z)}\right.}{\partial_x\,P_n(x,z)\left|_{x=x_{n,k}(z)}\right.}.
\end{equation}
On the other hand,  from  Proposition \ref{proL}
we get
\begin{equation}\label{diffP}
\partial_xP_{n}(x,z)\left |_{x=x_{n,k}(z)}\right.=\dfrac{P_{n-1}(x_{n,k}(z),z)}{A_{n-1}(x_{n,k}(z),z)},
\end{equation}
where $A_{n-1}(x,z)$ was defined in \eqref{operatorL}.
From the above we get
\begin{pro}
Let $(x_{n,k}(z))_{k=1}^n$ be the zeros of $P_{n}(x,z)$. Then they satisfy the differential equation
\begin{equation*}
\dot{x}_{n,k}(z)=\dfrac{x_{n,k}(z)}{z}\left[\dfrac{C_{n-1}(z,z)}{C_{n-1}(x_{n,k}(z),z)}\right],
\end{equation*}
where the polynomial $C_{n-1}(x,z)$ was defined in Proposition \ref{proL}.
\end{pro}
\begin{proof}
Taking into account \eqref{diffy}, \eqref{diffP}, \eqref{opediff2} and \eqref{operatorL}
$$\begin{aligned} \dot{x}_{n,k}(z)&=-A_{n-1}(x_{n,k}(z),z)\dfrac{\dot {P}_n(x,z)\left|_{x=x_{n,k}(z)}\right.}{P_{n-1}(x_{n,k}(z),z)}\\
&=\dfrac{A_{n-1}(x_{n,k}(z),z)}{z\,P_{n-1}(x_{n,k}(z),z)}\left(\dfrac{x_{n,k}(z)\,P_{n-1}(x_{n,k}(z),z)}{A_{n-1}((x_{n,k}(z),z))}-a_n(z)P_{n-1}(x_{n,k}(z),z)\right)\\
&=\dfrac{x_{n,k}(z)}{z}-a_{n}(z)\dfrac{A_{n-1}(x_{n,k}(z),z)}{z\,P_{n-1}(x_{n,k}(z),z)}\\
&=\dfrac{1}{z}\left(x_{n,k}(z)-\dfrac{(x_{n,k}(z)-z)\,x_{n,k}(z)}{x_{n,k}(z)+b_{n}(z)-(2n+\alpha+z+1)}\right)\\
&=\dfrac{x_{n,k}(z)}{z}\left[\dfrac{z+b_{n}(z)-(2n+\alpha+z+1)}{x_{n,k}(z)+b_{n}(z)-(2n+\alpha+z+1)}\right]
\end{aligned}$$
and the result follows.	
\end{proof}
\section{Concluding remarks}\label{section10}
We have studied the family of \textit{truncated Laguerre polynomials} $P_n(x; z)$, which are orthogonal with respect to the linear functional
$$\prodint{\elln,p}=\int_0^{z} p(x) x^\alpha e^{-x}dx, z>0, \quad  \alpha >-1.$$
Such a linear functional satisfies the Pearson equation
$$D((x-z)x\elln)+(-x(z+2+\alpha-x)+z(1-\alpha))\elln=0.$$ Using the symmetrization process, we construct the generalized truncated Hermite polynomials $S_n(x;z)$. These polynomials are semiclassical polynomials of class 2, if $\alpha=-1/2,$ and of class 3, if $\alpha\ne-1/2$. Taking into account the properties of truncated Laguerre and its closed relation with truncated Hermite  polynomials,  the Laguerre-Freud equations, the lowering and raising operators and the corresponding  holonomic differential equation that they satisfy are deduced. As an application, the equation of motion as well as the electrostatic interpretation of  the zeros of the above sequences of orthogonal polynomials are deduced.  We also obtained a second-order linear recurrence for the moments of $\elln$ and $\hn$, as well as a Painlev\'e differential equation with respect to the variable $z$ that the coefficients of the three term recurrence relation satisfy. Finally, differential equations  for the Stieltjes functions $\Su_{\elln}(t,z)$ and $\Su_{\hn}(t, z)$ are also deduced. Notice that the truncated Laguerre polynomials are related to Toda deformations of  the Jacobi weight $w(x)= x^{\alpha}, \alpha >-1$ in $(0,1).$

In a further work we are interested in asymptotic expansions for $P_n(x; z)$ and $S_n(x; z)$ as $n\to\infty$, $z\to \infty $ as well as when both $n, z \to \infty$ simultaneously. On the other hand, following Remark \ref{rem3,7} we will also focus our attention on the analysis of  generating functions for the truncated Laguerre and truncated generalized Hermite polynomials. Finally, the computation of zeros of $P_n(x; z)$ and their application in Gauss quadrature formulas following \cite{Gradimir18, Gradimir22} will be also done.

\section*{Acknowledgements}
The work of J. C. Garc\'ia-Ardila  has been supported by the Comunidad de Madrid multiannual agreement with the Universidad Rey Juan Carlos under the grant Proyectos I+D para Jóvenes Doctores, Ref. M2731, project NETA-MM. The work of F. Marcell\'an has been supported by FEDER/Ministerio de Ciencia e Innovación-Agencia Estatal de Investigación of Spain, grant PID2021-122154NB-I00, and the Madrid
Government (Comunidad de Madrid-Spain) under the Multiannual Agreement with UC3M in the line of Excellence of University Professors, grant EPUC3M23 in the context of the V PRICIT (Regional Program of Research and Technological Innovation).

\end{document}